\documentclass[onefignum,onetabnum]{siamart171218}

\usepackage[utf8]{inputenc}
\usepackage{lipsum}
\usepackage{amsfonts}
\usepackage{mathabx} 
\usepackage{graphicx}
\usepackage{epstopdf}
\usepackage{algorithmic}
\ifpdf
  \DeclareGraphicsExtensions{.eps,.pdf,.png,.jpg}
\else
  \DeclareGraphicsExtensions{.eps}
\fi

\newsiamremark{remark}{Remark}
\newsiamremark{hypothesis}{Hypothesis}
\crefname{hypothesis}{Hypothesis}{Hypotheses}
\newsiamthm{claim}{Claim}

\headers{Hanno Gottschalk and Marco Reese}{Hanno Gottschalk and Marco Reese}

\title{An Analytical Study in Multi Physics and Multi Criteria Shape Optimization}

\usepackage{amsopn}

\makeatletter
\newcommand*{\addFileDependency}[1]{
  \typeout{(#1)}
  \@addtofilelist{#1}
  \IfFileExists{#1}{}{\typeout{No file #1.}}
}
\makeatother

\newcommand*{\myexternaldocument}[1]{%
    \externaldocument{#1}%
    \addFileDependency{#1.tex}%
    \addFileDependency{#1.aux}%
}

\ifpdf
\hypersetup{
  pdftitle={An Analytical Study in Multi Physics and Multi Criteria Shape Optimization},
  pdfauthor={Hanno Gottschalk and Marco Reese}
}
\fi
\title{An Analytical Study in Multi Physics and Multi Criteria Shape Optimization}
\author{Hanno Gottschalk and Marco Reese}

\myexternaldocument{ex_supplement}

\usepackage[utf8]{inputenc}
\usepackage{color}
\usepackage{graphicx}
\usepackage{amsmath,amsfonts,mathabx}
\usepackage{setspace}

\DeclareMathOperator*{\argmin}{arg\,min}

\begin{document}

\setstretch{1.0}

\maketitle

\begin{abstract}
	A simple multi-physical system for the potential flow of a fluid through a shroud, in which a mechanical component, 
	like a turbine vane, is placed, is modeled mathematically. We then consider a multi criteria shape optimization problem,
	where the shape of the component is allowed to vary under a certain set of 2nd order H\"older continuous differentiable 
	transformations of a baseline shape with boundary of the same continuity class. As objective functions, we consider a 
	simple loss model for the fluid dynamical efficiency and the probability of failure of the component due to repeated 
	application of loads that stem from the fluid's static pressure. For this multi-physical system, it is shown that, under 
	certain conditions, the Pareto front is maximal in the sense that the Pareto front of the feasible set coincides with the Pareto 
	front of its closure. We also show that the set of all optimal forms with respect to scalarization techniques deforms 
	continuously (in the Hausdorff metric) with respect to preference parameters.
\end{abstract}



\section{Introduction}
The design of a mechanical component requires choosing a material and a shape. Often, a component serves a primary objective, 
but also requires a certain level of endurance. Material damage is caused by the loads that are imposed during service. 
The quest for an optimal design in the majority of cases therefore is at least a bi-criteria optimization problem and in many 
cases a multi criteria one \cite{ehrgott2005multicriteria}.  

In mechanical engineering, multi criteria optimization often comes along with coupled multi physics simulations. If we take the 
design of turbine blades as an example, the simulation of external flows and cooling air flows inside a blade have to be combined 
with a thermal and a mechanical simulation inside the blade \cite{sultanian2018gas}.   

Mathematical optimization is widely used in mechanical engineering; see, e.g., \cite{rao2019engineering,marler2004survey}. On the 
other hand, some directions of contemporary mathematical research -- like topology optimization (we refer, e.g., to 
\cite{Allaire2012shape,Bendsoe2003sigmund}) -- were initiated by mechanical engineers \cite{eschenauer1994bubble}.  
While the given field is interdisciplinary, from the mathematical point of view one would not only like to propose and analyse new 
optimization algorithms, but also to understand the existence and the properties of optimal solutions. While for mono-criteria 
optimization such a framework has been established \cite{Chenais75,Fujii,hasmaek03,Allaire2012shape,DelfZol11}, a general framework 
for multi criteria optimization is still missing; see however \cite{hasmaek03,doganay19,chir:mult:2018} for numerical studies 
addressing the topic.

Component life models from materials science are used to judge the mechanical integrity of a component after a certain number of load 
cycles; see, e.g.,  \cite{Werkstoffe}. These models often use deterministic life calculation which predicts failure at the 
point of the highest mechanical loading and thus involves the non-differentiable formation of a minimum life over all points on the 
component or component's surface, depending on whether we have a volume or a surface driven damage mechanism. In recent time, such models 
have been extended by probabilities of failure \cite{fedelich98,HertelVormwNotch2012,GottschSchmitz,DissSchmitz2014,LCFSchmid2014,ProbLCF2013,made2018combined,babuvska2019spatial}, 
which involve integrals over local functions of the stress tensor in the component or on the component's surface, 
respectively. This approach makes it possible to compute shape derivatives and gradients \cite{BucBut05,SokZol92,gottschalk2018adjoint,DissLaura} 
and therefore places component reliability in the context of shape optimization.
However, as remarked in \cite{GottschSchmitz}, the 
probability of failure as a objective functional requires more regular solutions as provided by the usual weak theory based on $H^1$ 
Sobolev spaces \cite{ErnGuerm04}. As we find here, this is also the case for simplistic models of fluid dynamical efficiency. As in 
previous works \cite{GottschSchmitz,BittGottsch,GBS_Ceramic2014}, we therefore apply a framework based on H\"older continuous classical solution 
spaces and extend it to multi criteria optimization.  

Within this general framework, we prove the existence of Pareto optimal designs; see also \cite{doganay19} for a related result in a 
different setting requiring less regularity. Here, however, we show how to use the graph compactness property \cite[Subsection 2.4]{hasmaek03}
along with the lower semi-continuity of all objective functionals to prove certain maximality properties of the non-dominated feasible 
points: Namely that the Pareto front in the set of feasible points \cite{ehrgott2005multicriteria} coincides with the Pareto front of 
the closure of the feasible points. Put in other words, each dominated design is also dominated by at least one Pareto optimal design.

We give a simplistic multi physical system as an example that fits the general framework. This mathematical model couples a potential 
flow with structural mechanics and is motivated from gas turbine engineering. We define two (rather singular) objective functionals, 
namely a aerodynamic loss based on the theory of boundary layers \cite{schlichting17} and furthermore, the probability of failure after a 
certain number of load cycles \cite{GottschSchmitz,made2018combined}. Each of these models includes nonlinear functions that depend on 
second derivatives of the solution after restriction to the boundary of the underlying PDE's domains. For this system, we prove that the 
assumptions of the general framework are fulfilled and we conclude that a maximal Pareto front exists in this case.

Multi criteria optimization relates to preferences of a decision maker \cite{ehrgott2005multicriteria,marler2004survey}. 
Here, we are interested in continuity properties of Pareto optimal shapes, when the preference is expressed by a parameter in a merit 
function, which, e.g., could be the weights in a weighted sum approach. The stability of the optimal solutions to such scalarization 
techniques in dependence of a parameter is already investigated in the literature; see, e.g., \cite{bank83,guddat89,karwat85,karwat87}, 
for finite dimensional and infinite dimensional spaces. Here, we show that our general framework is indeed suitable to prove certain continuity 
properties of the $\argmin$ sets of scalarized multi criteria optimization problems in the Hausdorff distance as a function of the 
scalarization -- or preference -- parameter. Such structural properties of the Pareto front for the first time are applied in the context of 
shape optimization.  

As this work is focused on the existence and the mathematical properties of the Pareto front in shape optimization, it 
does not contain any numerical and algorithmic contributions. However, numerical implementations for the computation of the reliability of 
failure can be found in \cite{schmitz2013risk,hanno17}. The adjoint method for computing the shape derivatives for such objective functionals 
has been applied in \cite{gottschalk2018adjoint,NumShapeCer2017,Gottschalk2019}. In \cite{doganay19,schultes2020hypervolume}, numerical studies 
of multi-criteria shape optimization are presented and Pareto fronts are explored. While the objective functionals in these two works fit 
into our general framework, they do not (yet) include aerodynamic losses (see however \cite{hanno20}). The results we present here on the continuity 
properties of Pareto optimal solutions under a change of scalarization can be seen as a theoretical underpinning of the Pareto tracing method 
presented in \cite{bolten20}.  

Our paper is organised as follows:
We introduce the physical systems, which underlie the multi criteria shape optimization problem we consider in Section 
\ref{section:multi-physics}. Afterwards, in Section \ref{section:multicrit}, we describe our framework for multi criteria shape  
optimization. By deriving uniform bounds for the solution spaces of the physical systems in Section \ref{section:existence}, we prove the 
well posedness of the shape optimization problem. Up to here we considered optimality in terms of Pareto optimality. In Section 
\ref{section:scalarization}, we apply scalarization techniques to transform the problem into an univariate shape problem and investigate the 
dependency of the optimal shapes on the specific used technique. In Section \ref{section:outlook}, we give a resume and an outlook on future 
research direction. Some technical details on H\"older functions and solutions of elliptic partial differential equations can be found in 
Appendix \ref{app}.

\section{\label{section:multi-physics}A Simple Multi Physics System}
We intend to optimize the shape of some component, e.g., a turbine vane, in terms of reliability and efficiency. Reliability depends on surface 
and volume forces acting on the component. In our setting, the component lies in a shroud and within the shroud a fluid is flowing past the 
component. Due to static pressure, the fluid imposes a surface force on the component. Hence, it is indispensable to include the fluid flow field
into the optimization process. At the same time, the component leads to frictional loss in the fluid that diminishes the efficiency of the design. 

In the following, we describe a simple model, which approximates the fluid flow in a simple way as potential flow and model frictional loss via a post 
processing step to the solution that is based on a simple model for the boundary layer. We also consider the effect of the fluid's mechanical loads 
to the component. As the static pressure takes the role of a boundary condition for the partial differential equation of linear elasticity, the 
internal stress fields of the component depend on the flow field, too. The component's fatigue life that results in the probability of failure, i.e.,
the formation of a fatigue crack, as a second objective functional.

\subsection{\bf{Potential Flow Equation}}
As component we consider a compact domain $\Omega \subset \mathbb{R}^3$ with $C^{k,\alpha}$ boundary -- where we always have that 
$k \in \mathbb{N}_{0}$ and $\alpha \in \, ]0,1]$ unless we specify further -- that is partially contained in some larger compact 
domain $D\subset\mathbb{R}^3$ representing a shroud with $C^{k,\alpha}$ boundary as well. With $\mathrm{int}(\Omega)$ and $\mathrm{cl}(\Omega)$
we shall denote the topological interior and closure of a set $\Omega$, respectively. 
We assume that
$\mathrm{int}(D\backslash\Omega)$ is simply 
connected, has $C^{k,\alpha}$ boundary, and that there exists an open ball $B := B_\epsilon$ with $\epsilon > 0$ such that 
$\mathrm{cl}(B) \subset \mathrm{int}(\Omega\backslash D)$. The shroud $D$ has an inlet and outlet where the fluid flows in and out, respectively. 
At the remaining boundary part the fluid cannot leak. In this work, we consider an incompressible and rotation free perfect fluid in a steady state. The 
assumption of zero shearing stresses in a perfect fluid -- or zero viscosity -- simplifies the equation of motion so that potential theory can be applied. 
The resulting solution still provides reasonable approximations to many actual flows. The viscous forces are limited to a thin layer of fluid adjacent to the 
surface, and therefore, in favor of simplicity, we leave these effects out since they have little effect on the general flow pattern\footnote{unless the local 
effects make the flow separate from the surface}.

A fundamental condition is that no fluid
can be created or destroyed within the shroud $D$. The equation of continuity expresses this condition. Consider a three dimensional
velocity field $v$ on $D\subset \mathbb{R}^3$, then the continuity equation is given by
\begin{equation*}
	\nabla\cdot v = 0.
\end{equation*}
If we assume that the velocity field $v$ is rotation free, $\nabla\times v=0$, then there exists a velocity potential or flow  potential $\phi$ such that
\begin{equation*}
	v = \nabla\phi.
\end{equation*}
Hence, under the assumption that $v$ is divergence free and rotation free, there is a velocity potential $\phi$ that
satisfies the Laplace equation
\begin{equation*}
	\Delta \phi =\nabla\cdot\nabla\phi= 0.
\end{equation*}
Let $n$ be the unitary outward normal of the boundary $\partial D$. By applying
suited Neumann boundary conditions $g$, that correspond with our assumptions for a conserved flow through the inlet and outlet of the shroud, we get 
the potential flow equation
\begin{equation}
\label{eq:pot-eq}
\begin{tabular}{ l l }
$\nabla v = \Delta \phi = 0$ & in $\mathrm{int}(D\backslash\Omega)$ \\
$v_n = v \cdot n = \frac{\partial \phi}{\partial n} = g$ & on $\partial D\backslash\partial(D \cap \Omega)$ \\
$v_n =  v \cdot n = \frac{\partial \phi}{\partial n} = 0$ & on $ D\cap\partial \Omega$. \\
\end{tabular}
\end{equation}
Here, we assume that $g$ is only non-zero in the inlet and outlet regions and is continued to be zero on the upper and lower wall of the shroud. 
Therefore, no discontinuities occur where $\partial\Omega$ meets $\partial D$. 

\begin{figure}[h]
	\centering
    \includegraphics[width=0.45\textwidth]{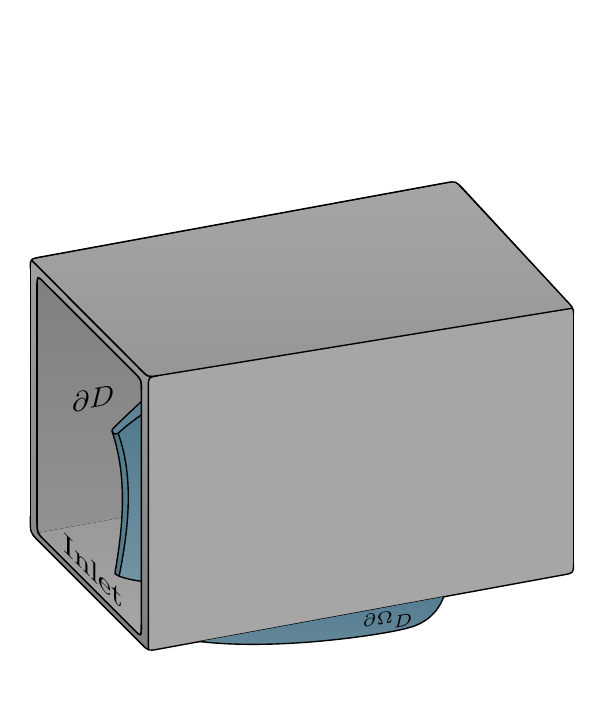}
    \includegraphics[width=0.45\textwidth]{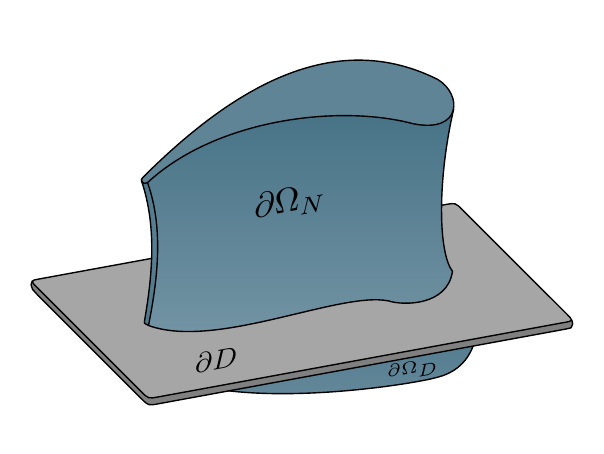}
	\caption{\setstretch{1.0} A turbine blade $\Omega$ within a shroud $D$. We note that this representation of the domains $\Omega$ and $D$ 
	is only a sketch. In particular, every visible edge has to be sufficiently rounded in order to ensure the Hölder continuity of the 
	boundaries.}
\end{figure}

The following lemma ensures the existence of a solution to the potential equation. It also gives a Schauder estimate
that leads to a uniform bound for the solution space we investigate in subsection \ref{subsection:uniform}. This uniform
bound is crucial for the existences of solutions to the multi criteria optimization problem we consider in this work and
which we introduce in Section \ref{section:multicrit}.

\begin{lemma}[Schauder Estimate for Flow Potentials]
	\label{lemma:pot-eq}
	We consider the potential flow equation {\normalfont(\ref{eq:pot-eq}).
	Let
	$g \in C^{1,\alpha}(D, \mathbb{R})$
	and assume that the boundaries described above are all of class $C^{k, \alpha}$ with $k \geq 2$. If $\int_{\partial D} g\, dA = 0$, then
	{\normalfont (\ref{eq:pot-eq})}
	possesses at least one solution 
	$\phi \in C^{2,\alpha}(\mathrm{cl}(D\backslash\Omega), \mathbb{R})$.
	To obtain uniqueness, we fix $u(x_{0}) = 0$ at some point $x_0 \in \mathrm{int}(D\backslash\Omega)$. 
	This solution satisfies
	\begin{equation}
	\label{est:pot-eq}
	\lVert \phi \rVert_{C^{2,\alpha}(D\backslash\Omega)} \leq
	C\left( \lVert \phi \rVert_{C^{0,\alpha}(D\backslash\Omega)} + 
	\lVert g \rVert_{C^{1,\alpha}(\partial D\backslash\partial(D \cap \Omega))} \right),
	\end{equation}
	with constant $C=C(\Omega)$.}
\end{lemma}
\begin{proof}
	By assumption (and by definition), $D\backslash\Omega$ has a Hölder continuous boundary of class $C^{k, \alpha}$, with 
	$k \geq 2,\, \alpha \in ]0,1]$. Additionally, the Neumann boundary condition holds $\int_{\partial(D \backslash \Omega)} g + 0 \, dA = 0$, and thus 
	the assertion follows directly out of \cite[Theroem 3.1 and Theorem 4.1]{nardi15}
\end{proof}

\subsection{\bf{Elasticity Equation}}
One of the most
crucial demands on the
component $\Omega$ is the reliability. Fatique failure is the most
appearing type of failure for, e.g., gas turbines where the event of failure
for a component as, e.g., a blade or vain is the appearance of the first crack.
For this purpose, we consider the elasticity equation, which models the
deformation of a component under given surface and volume forces and allows us to calculate the stress fields that drive crack formation. 

We denote with $n$ the unitary outward normal of the boundary $\partial\Omega$ and let
$\partial(\Omega \backslash B) = \partial\Omega \cup \partial B$ such that
$\partial B$ is clamped, and on $\partial\Omega$ a
force surface density $g_{\lvert \partial\Omega}$ is imposed. Then according
to \cite{ErnGuerm04} the mixed problem of linear isotropic elasticity, or the elasticity equation, is described by
\begin{equation}
\label{eq:el-eq}
\begin{tabular}{ l l }
$\nabla \cdot \sigma(u) + f = 0$ & in $\mathrm{int}(\Omega\backslash B)$ \\
$\sigma(u) = \lambda(\nabla \cdot u) I + \mu (\nabla u + \nabla u^T)$ & 
in $\mathrm{int}(\Omega \backslash B)$ \\
$u = 0$ & on $\partial B$ \\
$\sigma(u) n = g$ & on $\partial \Omega$ \\
\end{tabular}
\end{equation}
Here, $\lambda > 0$ and $\mu > 0$ are the Lam\'e constants (also called coefficients) and $u:\Omega\backslash B \to
\mathbb{R}^3$ is the displacement field on $\Omega \backslash B$. $I$ is the $3\times 3$
identity matrix. The linearized strain rate tensor $\epsilon(u):\Omega\backslash B \to
\mathbb{R}^{3\times 3}$ is defined as $\epsilon(u) = \frac{1}{2}(\nabla u + \nabla
u^\top)$. Approximate numerical solutions can be computed by a finite element
approach (see, e.g., \cite{ErnGuerm04} or \cite{HetEsl09}).

The potential equation (\ref{eq:pot-eq}) gives the velocity field at the part of the component's boundary $\partial\Omega$ 
that lies within the shroud $D$. Assuming that the total energy density, also denoted as stagnation 
pressure $p_{\text{st}}$, is constant at the inlet, we can derive the static pressure $p_\text{s}$ 
from Bernoulli's law
\begin{equation}
\label{eq:bernoulli}
p_{\text{st}} = \frac{1}{2}\rho \lvert \nabla \phi \rvert^2 + p_{\text{s}} \Leftrightarrow -p_{\text{s}} = 
\frac{1}{2}\rho \lvert \nabla \phi \rvert^2 - p_{\text{st}},
\end{equation}
where $\rho$ denotes the density of the fluid.
We consider the static pressure $p_{\text{s}}$ as surface load on the component $\Omega$.
Therefore, by continuously extending $p_{\text{s}}$ to be zero on $\partial \Omega\backslash D$, the surface load $g_{\text{s}}$ is given by
\begin{equation*}
	g_{\text{s}} = - p_{\text{s}}n = 
	\left(\frac{1}{2}\rho \lvert \nabla \phi \rvert^2 - p_{\text{st}}\right) n.
\end{equation*}
This yields as boundary condition on $\Omega$ for (\ref{eq:el-eq})
\begin{equation}
\label{el:g}
\sigma(u) n = g_{\text{s}} \quad \Leftrightarrow \quad
\sigma (u) n = \left(\frac{1}{2}\rho \lvert \nabla \phi \rvert^2 -
p_{\text{st}}\right) n.
\end{equation}
Hence, the displacement vector $u$ to the elasticity equation not only depends
on the shape $\Omega$ but also on the solution $\phi$ to the 
potential equation (\ref{eq:pot-eq}).

Assuming the boundary of $\Omega$ is of class $C^{k, \alpha}$, with $k \geq 1$ and $\alpha \in \, ]0,1]$, 
the unitary outward normal $n$ is a function in 
$C^{k-1, \alpha}(\partial \Omega, \mathbb{R}^{3})$. Further, let $\phi$ be in $C^{k, \alpha}(\partial \Omega, \mathbb{R})$.
Since we model an incompressible flow, the fluid density $\rho$ is constant as well as, by assumption, the stagnation pressure 
$p_{\mathrm{st}}$. Therefore, the static pressure $g_{\mathrm{s}}$ lies in $C^{k-1, \alpha}(\partial \Omega, \mathbb{R}^3)$ and 
the following lemma provides, for $k \geq 2$, existence and uniqueness of a solution $u$ along with a 
Schauder estimate.

\begin{lemma}[Schauder estimate for displacement fields, \cite{DissLaura}] 
	Consider the elasticity equation {\normalfont (\ref{eq:el-eq})}.
	Let $\Omega \subset \mathbb{R}^3$ be a compact domain that possesses
	$C^{k,\alpha}$-boundary. As volume load we consider $f \in
	C^{k,\alpha}(\Omega,\mathbb{R}^3)$ and as surface load $g_{\text{s}} \in
	C^{k+1,\alpha}(\partial \Omega,\mathbb{R}^3)$. 
	Then, the disjoint displacement-traction problem given by
	\begin{equation*}
	\begin{tabular}{ l l }
	$\nabla \cdot \sigma(u) + f = 0$ & in $\mathrm{int}(\Omega\backslash B)$ \\
	$\sigma(u) = \lambda(\nabla \cdot u) I + \mu (\nabla u + \nabla u^T)$ & 
	in $\mathrm{int}(\Omega\backslash B)$ \\
	$u = 0$ & on $\partial B$ \\
	$\sigma(u) n = g_{\text{s}}$ & on $\partial\Omega$ \\
	\end{tabular}
	\end{equation*}
	has a unique solution $u \in C^{k+2,\alpha}(\mathrm{cl}(\Omega \backslash B), 
	\mathbb{R}^3)$, which satisfies
	\begin{equation}
	\label{est:el-eq}
	\lVert u \rVert_{C^{k+2,\alpha}(\Omega \backslash B)} \leq
	C\left( \lVert f \rVert_{C^{k,\alpha}(\Omega)} + 
	\lVert g_{\text{s}} \rVert_{C^{k+1,\alpha}(\partial\Omega)} +
	\lVert u \rVert_{C^{0,\alpha}(\Omega \backslash B)} \right),
	\end{equation}
	with constant $C(\Omega)>0$.
\end{lemma} 

\subsection{\label{section:OptRelEff}\bf{Optimal Reliability and Efficiency}} 
Low cycle fatigue (LCF) driven surface crack initiation is particularly important for
the reliability of highly loaded engineering parts as turbine components \cite{schmitz2013risk,made2017probabilistic}. The design of 
such engineering parts therefore requires a model that is capable of accurately quantify risk levels for LCF crack initiation,
crack growth and ultimate failure. Here we refer to the model introduced in \cite{made2018combined} that models
the statistical size effect but also includes the notch support factor by using stress gradients arising from the 
elasticity equation (\ref{eq:el-eq}):
\begin{equation}
\label{eqa:J_prob}
	J_R(\Omega, u_\Omega):=
	\int_{\partial\Omega\cap D} \left( 
	\frac{1}{N_\text{det} (\nabla u_{\Omega},\nabla^2
		u_\Omega(x))}
	\right)^m dA.
\end{equation} 
$\Omega$ represents the shape of the component, $u_\Omega$ is the displacement field and the 
solution to the elasticity equation on $\Omega \backslash B$, $N_\text{det}$ is the deterministic number of
life cycles at each point of the surface of $\Omega$
and $m$ is the Weibull shape parameter. The proability of failure (PoF) after $t$ load cycles is then given as
$PoF(t)=1-e^{-t^m J_R(\Omega,u_\Omega)}$. Minimizing the probability of failure thus clearly is equivalent 
to minimizing $J_R(\Omega,u_\Omega)$.

For a detailed discussion including experimental validation we refer to \cite{hanno17}. We can 
apply this model as cost functional in order to optimize the component $\Omega$ with respect to
reliability.
 
Another primary objective of the component is the
efficiency that is connected with the viscosity
of the fluid flowing through the shroud. Viscosity is a measure
which describes the internal friction of a moving fluid.
In a laminar fluid the effect of viscosity is limited to a thin layer near
the surface of the component. The fluid does not slip along the surface, but
adheres to it. In the case of potential flow, there is a transition from zero
velocity at the surface to the full velocity which is present at a certain 
distance from the surface. The layer where this transition takes place
is called the boundary layer or frictional layer. The thickness of the 
boundary layer is not constant but (roughly)
proportional to the square root of the kinematic viscosity $\nu$ and is growing 
from the leading edge, the location where the fluid first impinge on the surface
of the component. 
Friction of the fluid on the surface leads to energy dissipation. A coefficient
for the inflicted local wall shear stress is given by
\begin{equation}
\label{coefficient:share-stress}
	\tau_w(x) = \frac{0.322\cdot\mu \lvert v\rvert^\frac{3}{2}}
	{\sqrt{\nu \cdot\text{dist}_{\text{LE}}(x)}},
\end{equation}
where we denote with $\lvert \cdot \rvert$ the Euclidean norm, $\mu$ is the viscosity, and $\text{dist}_{\text{LE}}$ the distance to
the leading edge along the component's surface $\partial\Omega$. For a detailed introduction 
to boundary layer theory one can see, e.g.,
\cite{schlichting17,boeswirth14}.
With this coefficient one can derive an estimate for the loss 
of power due to friction given by 
\begin{equation}
\label{functional:energy}
	J_E(\Omega, \phi_\Omega) := \int_{\partial\Omega\cap D}
	\lvert v_{\Omega} \rvert \tau_w\, dA.
\end{equation}

For the multi physics and multi criteria shape optimization problem we introduce in the next chapter, we realize above,
objective functionals that contain 
boundary integrals of second order derivatives of the solutions of second order elliptic BVPs. This can only be realised if one 
considers regular shapes and strong solutions. Further, we have to apply additional
assumptions on our shapes, i.e., a fixed and unique leading edge for all shapes in 
the shape space. This however, changes very little in our general analysis.

\section{A Multi Criteria Optimization Problem}
\label{section:multicrit}
In this section we introduce the multi criteria shape optimization problem, which is 
based on the above boundary value problems and cost functionals,
investigated in this work.  But before we introduce our multiphysics shape problem, we present a general approach 
to shape optimization problems and provide conditions for the existences of optimal shapes (see \cite[Section 2.4]{hasmaek03}). 
Afterwards, we apply our model.

\subsection{\bf{General Definitions}}
We denote a family of admissible shapes with
$\tilde{\mathcal{O}}$ and for every shape $\Omega\in \tilde{\mathcal{O}}$ we denote
with $V_1(\Omega),\dots ,V_n(\Omega),\, n\in\mathbb{N}$ state spaces of real valued functions on $\Omega$.
Consider a sequence of shapes $(\Omega_m)_{m\in \mathbb{N}}$ in $\tilde{\mathcal{O}}$, 
and let $\Omega \in \tilde{\mathcal{O}}$. Assuming a topology on the shape space $\tilde{\mathcal{O}}$ is given,
the convergence of $\Omega_m$ against $\Omega$ is denoted by 
$\Omega_m \overset{\tilde{\mathcal{O}}}{\longrightarrow} \Omega$ as $m \to \infty$. For a sequence of 
functions $(\boldsymbol{y}_m)_{m\in \mathbb{N}}$, with 
$\boldsymbol{y}_m \in 
\bigtimes_{i=1}^n V_i(\Omega_m)$ for all $m \in \mathbb{N}$, we 
denote the convergence against some $\boldsymbol{y} 
\in \bigtimes_{i=1}^n
V_i(\Omega)$ with $\boldsymbol{y}_m \rightsquigarrow 
\boldsymbol{y}$ as $m \to 
\infty$. We assume that for every $\Omega \in \tilde{\mathcal{O}}$ one can solve uniquely 
a given set of state problems, e.g., a set of PDEs or a variational inequalities. By associating the corresponding 
unique solutions $v_{i,\Omega}\in V_i(\Omega)$ with $\Omega \in \tilde{\mathcal{O}}$, one obtains the map 
$v_i: \Omega \mapsto v_{i,\Omega} \in V_i(\Omega)$. Let $\mathcal{O}$ be a subfamily of $\tilde{\mathcal{O}}$, then $\mathcal{G}
= \{ (\Omega,\boldsymbol{v}_\Omega)\, :\, \Omega \in \mathcal{O} \}$ is called the graph of the mapping $\boldsymbol{v} := (v_1,\dots, v_n)$. 
A cost functional $J$ on $\tilde{\mathcal{O}}$ is given by a map $J:(\Omega,\boldsymbol{y}) \mapsto 
J(\Omega,\boldsymbol{y}) 
\in \mathbb{R}$, where $\Omega \in \tilde{\mathcal{O}}$ and 
$\boldsymbol{y}\in \bigtimes_{i=1}^n V_i(\Omega)$. Then, a vector of
$l$ cost functionals is defined by $\boldsymbol{J}:=
(J_1,\dots, J_l)$, and the image of $\mathcal{O}$ (or $\mathcal{G}$) 
under $\mathbf{J}$ is denoted with $\mathcal{Y} \subset \mathbb{R}^l$.
For the sake of convenience, we shall write $\boldsymbol{J}(\Omega,\boldsymbol{v}_{\Omega})
:= (J_1(\Omega,\boldsymbol{v}_{\Omega}),\dots, J_l(\Omega,\boldsymbol{v}_{\Omega}))$, and, in addition, 
we make use of the notation $\nabla 
\boldsymbol{v}_{\Omega} := (\nabla v_{1,\Omega}, \dots, \nabla v_{n,\Omega})$.

\begin{definition}[Pareto optimality]
	Consider a subfamily $\mathcal{O}$ of $\tilde{\mathcal{O}}$
	with corresponding graph $\mathcal{G}$ to given state spaces
	$\boldsymbol{V} = (V_1, \dots,V_n)$.
	A point $(\Omega^*,\boldsymbol{v}_{\Omega^*}) \in \mathcal{G}$ is called
	Pareto optimal with respect to cost functionals
	$\boldsymbol{J} = (J_1,\dots, J_l)$,
	if there is no
	$(\Omega,\boldsymbol{v}_{\Omega})\in \mathcal{G}$ such that
	$J_k(\Omega,\boldsymbol{v}_{\Omega}) \leq J_k(\Omega^*,\boldsymbol{v}_{\Omega^*})$ for all
	$1\leq k \leq l$ and $J_i(\Omega,\boldsymbol{v}_{\Omega}) < 
	J_i(\Omega^*,\boldsymbol{v}_{\Omega^*})$ for some $i \in \{1,\dots,l\}$.
	The associated value $\mathbf{J}(\Omega^*,\boldsymbol{v}_{\Omega^*})$ is called
	nondominated.
\end{definition}
	
Let $\mathcal{Y}:= \mathbf{J}(\mathcal{G}) = 
\{ \mathbf{J}(\Omega, \boldsymbol{v}_{\Omega})\, :\, (\Omega, \boldsymbol{v}_{\Omega}) \in \mathcal{G} \}$ denote the image of the graph $\mathcal{G}$
under the objective functionals mapping $\mathbf{J}$.
For a set of Pareto optimal points,
we can define $\mathcal{Y}_N:= \{\mathbf{J}(\Omega, \boldsymbol{v}_{\Omega}) \in \mathcal{Y}\, :\, 
\mathbf{J}(\Omega, \boldsymbol{v}_{\Omega}) \text{ is nondominated in } \mathcal{Y}\}$, i.e.,
the corresponding Pareto front which lies by definition on the boundary of $\mathcal{Y}$.

\begin{definition}[Multi criteria shape optimization problem]
	Consider a subfamily $\mathcal{O}$ of $\tilde{\mathcal{O}}$ and for every $\Omega\in 
   	\mathcal{O}$ let $\boldsymbol{v}_\Omega = (v_{1,\Omega},\dots,v_{n,\Omega})$ be the unique solutions to given state problems on $\Omega$, 
	and let $\boldsymbol{J}= (J_1,\dots, J_l)$ be cost functionals on $\tilde{\mathcal{O}}$. 
	We define an optimal shape design problem by
	\begin{equation}
	\label{eq:opt_shape_design_problem}
	\begin{split}
	\left\{
	\begin{array}{ll}
	\text{Find } \Omega^* \in \mathcal{O} \text{ such that}\\
	(\Omega^*,\boldsymbol{v}_{\Omega^*})\text{ is Pareto optimal with respect to }
	\boldsymbol{J}.
	\end{array}
	\right.
	\end{split}
	\end{equation}
\end{definition}

The next theorem gives us conditions for the existence of a solution to the optimal
shape design problem (\ref{eq:opt_shape_design_problem}). Afterwards, in subsection 3.2, we define our shape optimization
problem and use this theorem to prove the existence of a solution to it.

\begin{theorem}
	\label{theorem:shape-solution-existence}
	Let $\tilde{\mathcal{O}}$ be a family of admissible domains and $\mathcal{O}$
	a subfamily. Consider cost functionals $\boldsymbol{J} =
	(J_1,\dots, J_l)$ on $\tilde{\mathcal{O}}$ and 
	assume for each $\Omega\in \tilde{\mathcal{O}}$ we have state problems with state
   	spaces $\boldsymbol{V}(\Omega) = (V_1(\Omega),\dots,V_n(\Omega))$ such that each state problem has a unique solution 
   	$v_{k,\Omega}\in V_k(\Omega)$, $1\leq k \leq n$. When the following both assumptions hold true
	\vspace*{1em}
	\begin{itemize}
		\item[(i)] Compactness of $\mathcal{G} = 
		\{ (\Omega,\boldsymbol{v}_\Omega)\, :\, \Omega \in \mathcal{O} \}$:\\
		Every sequence $(\Omega_m,\boldsymbol{v}_{\Omega_m})_{m\in\mathbb{N}}$ has a subsequence
		$(\Omega_{m_k}, \boldsymbol{v}_{\Omega_{m_k}})_{k\in\mathbb{N}}$ 
		that satisfies
		\begin{equation*}
			\begin{split}
				\Omega_{m_k} \overset{\tilde{\mathcal{O}}}{\longrightarrow} \Omega, \quad &k \to \infty \\
				\boldsymbol{v}_{\Omega_{m_k}} \rightsquigarrow 
				\boldsymbol{v}_\Omega, \quad &k\to \infty,
			\end{split}
		\end{equation*}
		for some $(\Omega,\boldsymbol{v}_\Omega) \in \mathcal{G}$.
		\vspace*{.5em}
		\item[(ii)] Lower semicontinuity of $J_k$:\\
		Let $(\Omega_m)_{m\in \mathbb{N}}$ be a sequence in $\tilde{\mathcal{O}}$ and
		$(\boldsymbol{y}_m)_{m\in \mathbb{N}}$ be a sequence such that $\boldsymbol{y}_m \in \boldsymbol{V}(\Omega_m)$ for all
		$m \in \mathbb{N}$. Consider some elements $\Omega,\,
		\boldsymbol{y}$ in
		$\tilde{\mathcal{O}}$ and $\boldsymbol{V}(\Omega)$, 
		respectively. Then,
		\begin{equation*}
			\left.
			\begin{array}{rl}
				\Omega_m \overset{\tilde{\mathcal{O}}}{\longrightarrow} \Omega, &\quad m\to \infty \\
				\boldsymbol{y}_m \rightsquigarrow 
				\boldsymbol{y}, &\quad m\to \infty
			\end{array}
			\right\}
			\Longrightarrow 
			\liminf_{n\to\infty}J_k(\Omega_n,\boldsymbol{y}_n) 
			\geq J_k(\Omega,\boldsymbol{y}),
		\end{equation*}
		for all $1\leq k \leq l$.
	\end{itemize}
	Then, the multi criteria shape design problem {\normalfont (\ref{eq:opt_shape_design_problem})} 
	possesses at least one solution and the Pareto front covers all nondominated 
	points in $\mathrm{cl}({\mathcal{Y}})$, i.e., $\mathcal{Y}_N = \mathrm{cl}({\mathcal{Y}})_N$, the set of non-dominated points in the closure of $\mathcal{Y}$.
\end{theorem}

\begin{proof}
	First, we prove the existence of an optimal shape. \cite[Theorem 2.10, p. 46]{hasmaek03} shows that, in this setting, a lower 
	semicontinuous cost functional possesses at least one minimal solution. 
	We apply this theorem, without loss of generality, to cost functional $J_1$ and minimize it
	on $\mathcal{G}$. Due to the compactness of $\mathcal{G}$ and the lower semicontinuity 
	of $J_1$, the resulting set of arguments of the minimum 
	$\argmin_{(\Omega,\boldsymbol{v}_\Omega)\in\mathcal{G}}J_1$ is also compact. Hence, we can
	again apply \cite[Theorem 2.10, p. 46]{hasmaek03} to the next cost functional $J_2$ and minimize it on
	$\argmin_{(\Omega,\boldsymbol{v}_\Omega)\in\mathcal{G}}J_1$. We continue
	this procedure until we minimized each cost functional on its preceding cost functionals 
	set of arguments of the minimum. The last set then contains at least one Pareto optimal point.
	
	For the second assertion, we recall that $\mathcal{Y}_N$ lies on the boundary of 
	$\mathcal{Y}$ and it follows directly that
	$\mathcal{Y}_N \subseteq \mathrm{cl}({\mathcal{Y}})_N$. Conversely, let 
	$\mathbf{J}(\Omega^*, \boldsymbol{v}_{\Omega^*}) \in \mathrm{cl}({\mathcal{Y}})_N$. Consider a sequence 
	$(\mathbf{J}(\Omega_n, \boldsymbol{v}_{\Omega_n}))_{n\in \mathbb{N}} \subset \mathcal{Y}$
	with $\mathbf{J}(\Omega_n, \boldsymbol{v}_{\Omega_n}) \to \mathbf{J}(\Omega^*, \boldsymbol{v}_{\Omega^*})$ 
	as $n\to \infty$. We assume that the corresponding sequence $(\Omega_{n}, \boldsymbol{v}_{\Omega_{n}})_{n \in\mathbb{N}} \subset \mathcal{G}$ converges to some
	$(\Omega, \boldsymbol{v}_{\Omega}) \in \mathcal{G}$ as well (since $\mathcal{G}$ is compact we can always find a
	subsequence). Due to the lower semicontinuity of $\boldsymbol{J}$, we have
	\begin{equation*}
		J_i(\Omega,\boldsymbol{v}_\Omega) \leq
		\lim_{n\to\infty}J_i(\Omega_n,\boldsymbol{v}_{\Omega_n}) = 
		J_i(\Omega^*, \boldsymbol{v}_{\Omega^*}), \quad \text{for all } 1\leq i \leq l.
	\end{equation*} 
	The Pareto optimality of $\mathbf{J}(\Omega^*, \boldsymbol{v}_{\Omega^*})$ gives that $\mathbf{J}(\Omega, \boldsymbol{v}_{\Omega}) = 
	\mathbf{J}(\Omega^*, \boldsymbol{v}_{\Omega^*})$,
	and since $\mathbf{J}(\Omega, \boldsymbol{v}_{\Omega}) \in \mathcal{Y}$,
	it follows that $\mathbf{J}(\Omega^*, \boldsymbol{v}_{\Omega^*}) \in \mathcal{Y}$ and therefore $\mathrm{cl}({\mathcal{Y}})_N 
	\subseteq \mathcal{Y}_N$.
	
\end{proof}

\subsection{\bf{Multi Physics Shape Optimization}}
\label{section:problem}
In the previous subsection we introduced a general framework of multi criteria shape optimization. We now state
a class of shape optimization problems that includes the multi physics shape optimization problem given by the 
coupled potential and elasticity equation as introduced in Section \ref{section:multi-physics}. As we will see, in Section 
\ref{section:existence}, multi criteria shape optimization problems from this class fulfill the required assumption of 
Theorem \ref{theorem:shape-solution-existence} to ensure us the existence of the Pareto front.

We consider shapes with Hölder continuous boundaries. 
This assumptions ensures, in this setting, strong regularity for the solutions of the physical problems
which enables us to deal with cost
functionals, defined on the boundaries of the shapes, containing first and second derivatives as motivated in subsection 
\ref{section:OptRelEff}. In the following, $C^{k,\alpha}$ stands for the real valued functions with $k$-th derivatives being 
H\"older continuous with exponent $\alpha$; see the appendix. \\

\begin{definition}
	Let $\Omega,\, \Omega^\prime$ be
	bounded domains in $\mathbb{R}^{d}$.
	\begin{itemize}
		\item[(i)] A $C^{k,\alpha}$-diffeomorphism from $\Omega$ to $\Omega^{\prime}$ is a bijective mapping
		$f: \Omega \to \Omega^\prime$ such that $f \in \left[
		C^{k,\alpha}(\Omega) \right]^d$ and $f^{-1} \in \left[
		C^{k,\alpha}(\Omega^\prime) \right]^d$.
		\item[(ii)] The set of $C^{k,\alpha}$-diffeomorphisms is denoted by
		$\mathcal{D}^{k,\alpha}(\Omega,\Omega^\prime)$ or
		$\mathcal{D}^{k,\alpha}(\Omega)$ if $f: \Omega \to \Omega$.
	\end{itemize}
\end{definition}

\begin{definition}
	\label{def:hemisphere}
	Consider a bounded domain $\Omega \subset \mathbb{R}^d$. The boundary of $\Omega$ is of class 
	$C^{k,\alpha}$, with $k\in \mathbb{N}_{0}$ and $0 < \alpha \leq 1$, if at each point $x_0 \in \partial \Omega$ there is an
	open ball $B=B(x_0)$ and a $C^{k,\alpha}$-diffeomorphism $T$ of $B$ onto $G \subset \mathbb{R}^d$ such that:
	\begin{equation*}
		(i)\ T(B \cap \Omega) \subset \mathbb{R}_{\geq 0}^d; \quad 
		(ii)\ T(B \cap \partial \Omega) \subset \partial \mathbb{R}^d_{\geq 0}; \quad
	\end{equation*}
	We shall say that the diffeomorphism $T$ straightens the boundary near $x_0$
	and call it hemisphere transform.
	Note that by this definition $\Omega$ is of class $C^{k,\alpha}$ if each point of 
	$\partial \Omega$ has a neighbourhood in which $\partial \Omega$ is the graph of a 
	$C^{k,\alpha}$ function of $d-1$ of the coordinates $x_1,\dots,x_n$. The converse
	is true if $k\geq 1$; see, e.g., \cite[Chapter 2, Theorem 5.5]{DelfZol11}. 
\end{definition}

\begin{definition}
	\label{def:shape-space}
	Let $K>0$ be a positive constant and $\Omega_0 \subset \Omega^{\text{ext}} \subset
	\mathbb{R}^{3}$ be compact $C^{k,\alpha}$ domains. The elements of the set
	\begin{align*}
		U_{k,\alpha}^{\text{ad}}(\Omega^{\text{ext}}) &:= \left\{
		\psi \in \mathcal{D}^{k,\alpha}(\Omega^{\text{ext}}) \, : \,
		\psi\lvert_{\mathrm{cl}(\Omega\backslash D)} \,= \text{id}, \,
		\lVert \psi \rVert_{\left[ C^{k,\alpha}(\Omega^{\text{ext}}) \right]^{3}}
		\leq K,\right.\\
		&\hspace{6cm}\left.\lVert \psi^{-1} \rVert_{\left[ C^{k,\alpha}(\Omega^{\text{ext}}) \right]^{3}}
		\leq K
		\right\}
	\end{align*}
	are called design-variables. These design variables induce, in a natural way,
	the set of admissible shapes
	\begin{equation*}
		\mathcal{O}_{k,\alpha} := \mathcal{O}_{k,\alpha}(\Omega_0,\Omega^{\text{ext}}) :=
		\left\{ \psi(\Omega_0)\, : \, \psi \in U_{k,\alpha}^{\text{ad}}(\Omega^{\text{ext}})
		\right\},
	\end{equation*}
	assigned to $\Omega_0$. Note that due to the Hölder continuity, every $\Omega\in
	\mathcal{O}_{k,\alpha}$ is compact.
\end{definition}

\begin{lemma}\label{lemma:cone_condition}
	Let $k \geq 0$ and $\alpha \in \, ]0,1]$, then the shape space $\mathcal{O}_{k, \alpha}$ satisfies a 
	uniform cone condition.
\end{lemma}
\begin{proof}
	As $k \geq 1$, the shape $\Omega_0$ is a domain with Lipschitz boundary
	and therefore fulfills
	a uniform cone condition. Since every transform $\psi \in U_{k,\alpha}^{\text{ad}}(\Omega^{\text{ext}})$ 
	is a $C^{k, \alpha}$-diffeomorphism with 
	$\lVert \psi \rVert_{\left[ C^{k,\alpha}(\Omega^{\text{ext}}) \right]^{3}}\leq K$
	and $\lVert \psi^{-1} \rVert_{\left[ C^{k,\alpha}(\Omega^{\text{ext}}) \right]^{3}}\leq K$, we have
	\begin{equation}\label{eq:cone_equation}
		\frac{1}{K} |x - y| \leq |\psi(x) - \psi(y)| \leq K |x - y|, \quad \text{for all } x,y \in \Omega_0.
	\end{equation}
	Let $C(x)$ be the cone associated with the cone condition satisfied by $\Omega_0$, where $x \in \partial \Omega_0$ 
	denotes the vertex. Further, we denote with $C_{K}(x)$ the cone where we decreased the radius of $C$
	with factor $\frac{1}{K}$. Then, by the lower bound in (\ref{eq:cone_equation}), we can always place the shrinked 
	cone $C_{K}$ within the transformed cone $\psi(C(x))$ at the boundary point $\psi(x)$ for every 
	$\psi \in U_{k,\alpha}^{\text{ad}}(\Omega^{\text{ext}})$ and $x \in \partial \Omega_0$. 
	Therefore, the cone $C_K$ provides the uniform cone condition for $\mathcal{O}_{k, \alpha}$.
	
\end{proof}

On the space of admissible domains we can define the Hausdorff distance as
metric. Note that in general the Hausdorff distance is no metric because in general
the identity of indiscernibles is not given. Here however, the compactness of
the shapes ensure us this property.

\begin{definition}
	For two non-empty subsets $\Omega$, $\Omega^\prime$ of a metric space $(M, d)$ 
	we define their Hausdorff distance by
	\begin{equation*}
		d_H (\Omega,\Omega^\prime) :=
		\max \{ \sup_{x\in \Omega} \inf_{y\in \Omega^\prime} d(x - y),
		\, \sup_{y\in \Omega^\prime} \inf_{x\in \Omega} d( x - y ) \}
	\end{equation*}
\end{definition}
If we equip the set $F(M)$ of all closed subsets of a metric space $(M,d)$ with the Hausdorff distance,
then we obtain another metric 
space. Since the shapes in $\mathcal{O}_{k,\alpha}$ are compact, the Hausdorff distance defines
a metric on $\mathcal{O}_{k,\alpha}$. By the following Lemma, we see in chapter 4, that
$(\mathcal{O}_{k,\alpha},\, d_H)$ is additionally compact.

\begin{theorem}[Blaschke's Selection Theorem \cite{baley40}]
	\label{thm:blaschke}
	Let $(M,d)$ be a metric space, where $M$ is a compact subset of a Banach space $B$. Then,
	the set $F(M)$ of all closed subsets of $M$ is compact with respect to the Hausdorff distance 
	$d_H$. 
\end{theorem}

\begin{definition}[Local Cost Functionals]
	\label{def:cost-functional}
	Let $\mathcal{O}\subset \mathcal{P}(\mathbb{R}^3)$
	denote a shape space with corresponding state spaces $V_1(\Omega),\dots, V_n(\Omega),
	\, \Omega \in \mathcal{O}$ and graph 
	$\mathcal{G} := \{(\Omega,\boldsymbol{v}_{\Omega})\,:\, \Omega\in \mathcal{O}\}$. Assuming that $V_i(\Omega) \subseteq 
	C^k(\Omega,\mathbb{R}^3)$ for all $1\leq i\leq n$, 
	the local cost functional on $\mathcal{G}$ is given by
	\begin{equation}
	\begin{split}
	J(\Omega,\boldsymbol{v}_{\Omega}) := &\int_\Omega \mathcal{F}_{\text{vol}}\hspace{.1em} 
	(x,\boldsymbol{v}_{\Omega},\nabla \boldsymbol{v}_{\Omega},\dots, \nabla^k \boldsymbol{v}_{\Omega})\,dx \\
	&+ \int_{\partial \Omega} 
	\mathcal{F}_{\text{sur}}\hspace{.1em}
	(x,\boldsymbol{v}_{\Omega},\nabla \boldsymbol{v}_{\Omega},\dots, \nabla^k \boldsymbol{v}_{\Omega})\,dA, 
	\end{split}
	\end{equation}
	where $\mathcal{F}_{\text{vol}}$,
	$\mathcal{F}_{\text{vol}}: \mathbb{R}^d\to \mathrm{cl}({\mathbb{R}}_{\geq 0})$ and $d =
	3 + n\sum_{j=0}^k 3^{j+1} = 3+\frac{3n}{2}(3^{k+1}-1)$. We 
	denote the volume integral and surface integral with
	\begin{equation*}
		\begin{split}
			&J_{\text{vol}}(\Omega,\boldsymbol{v}_{\Omega}) := \int_\Omega 
			\mathcal{F}_{\text{vol}}\hspace{.1em} 
			(x,\boldsymbol{v}_{\Omega},\nabla \boldsymbol{v}_{\Omega},\dots, \nabla^k \boldsymbol{v}_{\Omega})\,dx, \\
			&J_{\text{sur}}(\Omega,\boldsymbol{v}_{\Omega}) := \int_{\partial \Omega} 
			\mathcal{F}_{\text{sur}}\hspace{.1em}
			(x,\boldsymbol{v}_{\Omega},\nabla \boldsymbol{v}_{\Omega},\dots, \nabla^k \boldsymbol{v}_{\Omega})\,dA.
		\end{split}
	\end{equation*}
\end{definition}

\begin{definition}[Multi Physics Shape Optimization Problem]
	\label{def:multi-crit-problem}
	We consider the space
	$(\mathcal{O}_{k,\alpha}, d_H)$ of admissible shapes
	and let $J_1, \dots, J_l$ be local cost functionals on the Graph 
	$\mathcal{G} := \{
	(\Omega, u_\Omega, \phi_\Omega) \, : \, \Omega \in \mathcal{O}_{k,\alpha},\, u_\Omega
	\text{ solves } (\ref{eq:el-eq})\text{ on } \Omega, \\ \phi_\Omega \text{ solves }
	(\ref{eq:pot-eq}) \text{ on } \Omega \backslash D \}$.
	The multi physics shape optimization problem is given by:
	
	\begin{equation}
	\label{eq:multiphysics-shape-eq}
	\begin{split}
	\left\{
	\begin{array}{ll}
	\text{Find } \Omega^* \in \mathcal{O}_{k,\alpha} \text{ such that}\\
	(\Omega^*,u_{\Omega^*},\phi_{\Omega^*})\text{ is Pareto optimal with respect to }
	\boldsymbol{J}.
	\end{array}
	\right.
	\end{split}
	\end{equation}
\end{definition}

Before this section ends, we note that our choice of state 
problems here is only exemplary. In the next section we see
that we can include an arbitrary amount of physical models
in this multi physics shape optimization problem as long as they
provide an unique solution with sufficient regularity and a compact
solution space on the shape space $\mathcal{O}_{k,\alpha}$.

\section{\label{section:existence}Existence of Pareto Optimal Shapes}
In this section, the approach outlined in Thoerem \ref{theorem:shape-solution-existence} will be followed in order
to show the existence of an optimal shape for the multi physics shape optimization problem.
This approach includes the compactness of the Graph $\mathcal{G}$ which requires
bounded solution spaces on $\mathcal{O}_{k,\alpha}$. We first derive
such uniform bounds based on the Schauder estimates given in Section \ref{section:multi-physics}.
\subsection{\bf{Uniform Bounds for Solution Spaces}}

\label{subsection:uniform}

\begin{lemma}
	\label{lemma:pot-uni}
	Let $\phi_{\Omega}$ be the unique solution to {\normalfont (\ref{eq:pot-eq}) } on $\Omega \in \mathcal{O}_{k,\alpha}$ with $k\geq2$. 
	Then, there exists a constant $K>0$ independent of $\Omega$ such that
	\begin{equation*}
		\lVert \phi_{\Omega} \rVert_{C^{2,\alpha}(D\backslash\Omega)} \leq K.
	\end{equation*}
\end{lemma}

\begin{proof}
	This estimate is based on Lemma \ref{lemma:pot-eq}:
	\begin{equation*}
	\lVert \phi_{\Omega} \rVert_{C^{2,\alpha}(D\backslash\Omega)} \leq
	\tilde{C}\left( \lVert \phi_{\Omega} \rVert_{C^{0,\alpha}(D\backslash\Omega)} + 
	\lVert g \rVert_{C^{1,\alpha}(\partial D\backslash\partial(D \cap \Omega)} \right),
	\end{equation*}
	where the constant $\tilde{C}$ possibly depends on the shape $\Omega$. 
	First, we outline that the constant 
	$\tilde{C}$ can be chosen independently of the shape $\Omega\in\mathcal{O}_{k,\alpha}$. 
	A full proof is provided in \cite{GilbTrud} 
	(or \cite{BittGottsch}). 
	In order to prove estimate (\ref{est:pot-eq}), one straightens
	the boundary $\partial\Omega$ piecewise with hemisphere transforms. The dependence of the constant $\tilde{C}$ is through
	the ellipticity of the differential operator and hence 
	depends on the bounds of the hemisphere transform 
	that is used to straightens the boundary. Let $T$ be such
	a hemisphere transform for $\Omega_0$.
	For every shape
	$\psi(\Omega_0)\in\mathcal{O}_{k,\alpha}$ we can 
	construct a hemisphere transform by pulling $\psi(\Omega_0)$
	back to $\Omega_0$ and apply
	$T$ afterwards, i.e., $T_{\psi(\Omega_0)} := T \circ \psi^{-1}$.  
	Due to the definition of design variables, $\psi$
	is uniformly
	bounded w.r.t. $\lVert \cdot \rVert_{C^{k,\alpha}(\Omega^\text{ext})}$, and since $k \geq 2$, also
	lies in $C^{k-1, 1}(\Omega^{\text{ext}})$. Therefore, it follows 
	that $T_{\psi(\Omega_0)}$ is a function of class $C^{k, \alpha}$ and uniformly bounded in $\mathcal{O}_{k,\alpha}$ 
	w.r.t. $\lVert \cdot \rVert_{C^{k, \alpha}(\Omega^{\text{ext}})}$ as well.
	
	Next, we note that $\lVert g \rVert_{C^{1,\alpha}(\partial D\backslash\partial(D \cap \Omega)}$ is 
	obviously bounded by 
	$\lVert g \rVert_{C^{1,\alpha}(\partial D)}$,
	and it remains to further estimate $\lVert \phi_{\Omega} \rVert_{C^{0,\alpha}(D\backslash\Omega)}$.
	Since $\mathcal{O}_{k, \alpha}$ satisfies a uniform cone condition (see Lemma \ref{lemma:cone_condition}),
	\cite[Lemma 5.5]{GottschSchmitz} implies that for every $\epsilon > 0 $ there is a constant $C(\epsilon) > 0$
	such that
	\begin{equation*}
		\lVert \phi_{\Omega} \rVert_{C^{0,\alpha}(D\backslash\Omega)} \leq 
		\epsilon \lVert \phi_{\Omega} \rVert_{C^{1,\alpha}(D\backslash\Omega)} + 
		C(\epsilon) \int_{D\backslash\Omega} \lvert \phi_{\Omega} \rvert\, dx.
	\end{equation*}
	We choose $\epsilon < 1/\tilde{C}$ and get 
	\begin{equation*}
		\begin{split}
			\lVert \phi_{\Omega} \rVert_{C^{2,\alpha}(D\backslash\Omega)} \leq \
			&\tilde{C}\left( \lVert \phi_{\Omega} \rVert_{C^{0,\alpha}(D\backslash\Omega)} + 
			\lVert g \rVert_{C^{1,\alpha}(\partial D)} \right),\\
			\leq \ &\frac{1}{1 - \epsilon \tilde{C}} \left( C(\epsilon) \int_{D\backslash\Omega} \phi_{\Omega}\, dx +
			\tilde{C} \lVert g \rVert_{C^{1,\alpha}(\partial D)} \right) \\
			\leq \ &\frac{1}{1 - \epsilon \tilde{C}} \left( C(\epsilon) 
			\lVert \phi_{\Omega} \rVert_{H^1(D\backslash\Omega)} +
			\tilde{C} \lVert g \rVert_{C^{1,\alpha}(\partial D)} \right).
		\end{split}
	\end{equation*}
	One can easily verify the a-priori estimate $\lVert \phi_{\Omega} 
	\rVert_{H^1(D\backslash\Omega)} \leq C_p \lvert D \rvert
	\lVert g \rVert_{C^{1,\alpha}(\partial D)}$ 
	holds for a constant
	$C_p >0$ originating from the Poincar\'e inequality (see, e.g., \cite[Lemma B.61]{ErnGuerm04}). This yields
	\begin{equation*}
		\begin{split}
			\lVert \phi_{\Omega} \rVert_{C^{2,\alpha}(D\backslash\Omega)} \leq \
			&\frac{1}{1 - \epsilon \tilde{C}} \left( 
			C(\epsilon) \lVert \phi_{\Omega} \rVert_{H^1(D\backslash\Omega)} +
			\tilde{C} \lVert g \rVert_{C^{1,\alpha}(\partial D)} \right)\\ 
			\leq \ &\frac{C(\epsilon)C_p\lvert D \rvert + \tilde{C}}{1 - \epsilon \tilde{C}}
			\lVert g \rVert_{C^{1,\alpha}(\partial D)}
			=: K.
		\end{split}
	\end{equation*}
	\
\end{proof}

We recall that the elasticity equation (\ref{eq:el-eq}) 
describes with $g_{\text{s}}$, given by the static pressure
$p_{\text{s}}$, the surface force that the fluid exerts on the component $\Omega$ (see (\ref{el:g})).
Hence, in our framework, the surface force $g_{\text{s}}$ is given by Bernoulli's equation and we have
\begin{equation}
\sigma(u_{\Omega})n = g_{\text{s}} \quad \Leftrightarrow \quad
\sigma (u_{\Omega}) n = \left(\frac{1}{2}\rho \lvert \nabla \phi_{\Omega} \rvert^2 -
p_{\text{st}}\right) n,
\end{equation}
where we contiuously extends $g_{\text{s}}$ to be zero on $\partial \Omega \backslash D$.
The solution $u_{\Omega}$ of the elasticity equation not only depends on the 
shape $\Omega \in \mathcal{O}_{k,\alpha}$ but on the solution $\phi_{\Omega}$ of potential equation
(\ref{eq:pot-eq}) as well. We can derive a uniform bound for $u_{\Omega}$, as we have for 
$\phi_{\Omega}$, from estimate (\ref{est:el-eq}) which already provides an 
uniform bound in $\mathcal{O}_{k,\alpha}$
if the surface load $g_{\text{s}}$ is independent of $\phi_{\Omega}$. However, in
our framework $g_{\text{s}}$ depends 
on $\phi_{\Omega}$, and thus we have to further estimate the surface force $g_{\text{s}}$ in $\mathcal{O}_{k,\alpha}$. 

\begin{lemma}
	Let $u_{\Omega}$ be the unique solution to {\normalfont (\ref{eq:el-eq})} on $\Omega \in \mathcal{O}_{k,\alpha}$ with $k \geq 2$. 
	Then, there exists a constant $K>0$ independent of $\Omega$ such that 
	\begin{equation*}
		\lVert u \rVert_{C^{2,\alpha}(\Omega \backslash B)} \leq K,
	\end{equation*}
\end{lemma}

\begin{proof}
	We consider estimate (\ref{est:el-eq}) given by
	\begin{equation*}
	\lVert u_{\Omega} \rVert_{C^{2,\alpha}(\Omega\backslash B)} \leq
	\tilde{C}\left( \lVert f \rVert_{C^{0,\alpha}(\Omega)} + 
	\lVert g_{\text{s}} \rVert_{C^{1,\alpha}(\partial\Omega)} +
	\lVert u_{\Omega} \rVert_{C^{0,\alpha}(\Omega\backslash B)} 
	\right),
	\end{equation*}
	where the constant $\tilde{C}>0$ potentially depends on the domain $\Omega$.
	However, as we described above, due to the construciton of $\mathcal{O}_{k,\alpha}$, $C$ can be choosen 
	independently of $\Omega \in \mathcal{O}_{k,\alpha}$.
	Further, $\lVert f \rVert_{C^{0,\alpha}(\Omega)}$ is bounded by
	$\lVert f \rVert_{C^{0,\alpha}(\Omega^{\text{ext}})}$, and $\lVert g_{\text{s}} \rVert_{C^{1,\alpha}(\partial\Omega)}$ 
	depends on the potential $\phi_{\Omega}$ in terms of (\ref{el:g}). In the proof to Lemma \ref{lemma:pot-uni},
	we have shown for the bounded domains $\Omega \in \mathcal{O}_{k,\alpha}$ that 
	the diffeomorphisms, which describe the boundary of the domains, are uniformly bounded in $\mathcal{O}_{k,\alpha}$ with respect 
	to $\lVert \cdot \rVert_{C^{k,\alpha}}$.
	As we can use these diffeomorphism as chart mapping to describe the two-dimensional submanifold $\partial \Omega$, we can 
	conclude that the unitary normal vector $n$ of $\partial \Omega$ is uniformly bounded in $\mathcal{O}_{k,\alpha}$ by some 
	constant $M>0$ with respect to $\lVert \cdot \rVert_{C^{k, \alpha}(\partial \Omega)}$. Hence, we can estimate
	\begin{equation*}
		\begin{split}
			\lVert g_{\text{s}} \rVert_{C^{1,\alpha}(\partial\Omega)} = \
			&\left\lVert \left(\frac{1}{2}\rho \lvert \nabla \phi_{\Omega} \rvert^2 - 
			p_{\text{st}}\right) n \right\rVert_{C^{1,\alpha}(\partial\Omega)} \\
			\leq  \ &p_{\text{st}}M + \frac{1}{2}\rho  \lVert \nabla \phi_\Omega^2 n \rVert_{C^{1, \alpha}(\partial \Omega)}.
		\end{split}
	\end{equation*}
	Equation (\ref{eq:pot-eq}) models an incompressible fluid and therefore the fluid
	density $\rho$ is constant as well as $p_{\text{st}}$ by assumption. Since $n$ is uniformly bounded in $\mathcal{O}_{k, \alpha}$,
	we have by Lemma \ref{lemma:pot-uni} that $\lVert \nabla \phi_\Omega^2 n \rVert_{C^{1, \alpha}(\partial \Omega)}$ is also 
	uniformly bounded in $\mathcal{O}_{k, \alpha}$. Thus, we get that
	\begin{equation*}
		\lVert g_{\text{s}} \rVert_{C^{1,\alpha}(\partial\Omega)} \leq \tilde{M},
	\end{equation*}
	with some constant $\tilde{M}>0$ which is independent of $\Omega \in \mathcal{O}_{k,\alpha}$.
	
	Now, by \cite[Lemma 5.5]{GottschSchmitz}, for $\epsilon > 0$ one can estimate
	\begin{equation*}
		\lVert u_{\Omega} \rVert_{C^{0,\alpha}(\Omega\backslash B)} \leq 
		\epsilon \lVert u_{\Omega} \rVert_{C^{1,\alpha}(\Omega\backslash B)} + 
		C(\epsilon) \int_{\Omega\backslash B} \lvert u_{\Omega} \rvert\, dx,
	\end{equation*}
	with constant $C(\epsilon) > 0$ which is independent of $\Omega \in \mathcal{O}_{k,\alpha}$. Applying this with
	$\epsilon < 1/\tilde{C}$ on (\ref{est:el-eq})
	and estimating $\int_{\Omega\backslash B} \lvert u_{\Omega} \rvert\, dx \leq 
	\Vert u_{\Omega} \rVert_{H^1(\Omega\backslash B)}$
	yields
	\begin{equation*}
		\begin{split}
			\lVert u_{\Omega} \rVert_{C^{2,\alpha}(\Omega\backslash B)} \leq
			\frac{1 + \tilde{C}}{1-\epsilon\tilde{C}}\left( \lVert f \rVert_{C^{0,\alpha}(\Omega^{\text{ext}})} + 
			\lVert g_{\text{s}} \rVert_{C^{1,\alpha}(\partial\Omega)} +
			C(\epsilon)\lVert u_{\Omega} \rVert_{H^1(\Omega\backslash B)}\right)
		\end{split}
	\end{equation*} 
	Let $V_{DN} = \{ v \in [H^1(\Omega\backslash B)]^3\, :\, v=0 \text{ a.e.  on } 
	\partial B\}$, and consider the weak formulation of (\ref{eq:el-eq}) given by
	\begin{equation*}
		\int_{\Omega\backslash B} \text{tr}(\sigma(u_{\Omega})\epsilon(v))\, dx = 
		\int_{\Omega\backslash B} fv\, dx + \int_{\partial\Omega}g_{\text{s}}v\, dA, \quad \text{for all } v \in V_{DN}.
	\end{equation*}
	One can see that for all $v\in V_{DN}$ we have
	\begin{equation*}
		C^{-1}\lVert \epsilon(v) \rVert^2_{L_2(\Omega \backslash B)}\leq \int_{\Omega\backslash B} \text{tr}(\sigma(v)\epsilon(v))\, dx,
	\end{equation*}
	with constant $C>0$. We can choose $C$ such that, due to the uniform boundedness of $f$ and $g_{\text{s}}$,
	we also have
	\begin{equation*}
		\left\lvert \int_{\Omega\backslash B} fv\, dx + \int_{\partial\Omega}g_{\text{s}}v\, dA \right\rvert
		\leq C\lVert v \rVert_{H^1(\Omega\backslash B)},
	\end{equation*}
	where $C$ is additionally uniform with respect to $\mathcal{O}_{k,\alpha}$. 
	Korn's second inequality (\ref{lemma:Korn-2nd}) then implies
	\begin{equation*}
		\begin{split}
			q\lVert \epsilon(u_{\Omega}) \rVert^2_{L_2(\Omega\backslash B)}  \leq \,
			&C\lVert u_{\Omega} \rVert_{H^1(\Omega\backslash B)} \leq \,
			C_K \lVert \epsilon(u_{\Omega}) \rVert_{L_2(\Omega\backslash B)} \\
			\Rightarrow  \, &\lVert u_{\Omega} \rVert_{H^1(\Omega\backslash B)} \leq \,C_K^2,
		\end{split}
	\end{equation*}
	where the constant $C_{K} > 0$, which originates from Korn's second inequality, may depends on the domain $\Omega$.
	Examining the proof to Korn's second inequality (see, e.g., \cite{Nitsche1981}),
	one can see that the constant $C_K$ depends on $\Omega$ through the cone associated to the uniform cone condition to $\Omega$. 
	Since $\mathcal{O}_{k, \alpha}$ satisfies a uniform cone condition, the constant $C_K$ is uniform with respect to $\mathcal{O}_{k, \alpha}$.
	Thus, the previous inequality is uniform in $\mathcal{O}_{k,\alpha}$ and the assertion is proven.
	
\end{proof}

\subsection{\bf{Pareto Optimality}}
In order to prove the existence of an optimal shape to the multi 
physics shape optimization problem (\ref{eq:multiphysics-shape-eq}),
we want to make use of Theorem \ref{theorem:shape-solution-existence}. Therefore, we show that the local cost functionals from
Definition \ref{def:cost-functional} are lower semicontinuous --we even show that they are continuous--
and that the graph from Definition \ref{def:multi-crit-problem} is compact. The continuity is given
and discussed in Lemma \ref{lemma:cost-funct-continuity}. First, we denote with 
$\mathcal{P}_{k,\alpha}:= \{ \phi_\Omega \,:\, \phi_{\Omega} 
\text{ solves } (\ref{eq:pot-eq}) \text{ with } \Omega \in \mathcal{O}_{k,\alpha} \}$ and
$\mathcal{E}_{k,\alpha}:= \{ u_\Omega \,:\, u_{\Omega} 
\text{ solves } (\ref{eq:el-eq}) \text{ with } \Omega \in \mathcal{O}_{k,\alpha} \}$
the spaces of solutions to (\ref{eq:pot-eq}) and (\ref{eq:el-eq}) on admissible shapes, 
respectively. We equip these spaces
with the metric that is induced by the Hölder norm.
The solutions in these spaces are defined on different and distinct domains and therefore are
not comparable with respect to $\lVert \cdot \rVert_{C^{k,\alpha}}$. We give our solution to this
problem in the following first definition of this subsection:

\begin{definition}[$C^{k,\alpha}$-Convergence of Functions with Varying Domains]
\label{def:extensions}
	We recall the sets $\mathcal{O}_{k,\alpha}$ and $\Omega^{\text{ext}}$ from 
	Definition \ref{def:shape-space}. With $p_\Omega : [C^{k,\alpha}(\Omega\backslash B)]^{3n} \to 
	[C^{k,\alpha}_0(\Omega^{\text{ext}}\backslash B)]^{3n}$ we denote the extension operator that can be derived
	from Lemma \ref{lemma:extension}. For $\boldsymbol{v}
	\in [C^{k,\alpha}(\Omega\backslash B)]^{3n}$ set $\boldsymbol{v}^{\text{ext}} 
	= p_\Omega \boldsymbol{v}$. For $(\Omega_n)_{n\in\mathbb{N}} \subset \mathcal{O}_{k,\alpha},\ 
	\Omega \in \mathcal{O}_{k,\alpha}$ and $(\boldsymbol{v}_n)_{n\in \mathbb{N}}$ with 
	$\boldsymbol{v}_n \in 
	[C^{k,\alpha}(\Omega_n\backslash B)]^{3n},\ n\in \mathbb{N}$, the expression $\boldsymbol{v}_n \rightsquigarrow
	\boldsymbol{v}$ as $n \to \infty$ is defined 
	by $\boldsymbol{v}_n^{\text{ext}} \to \boldsymbol{v}^{\text{ext}}$ in 
	$[C^{k,\alpha}_0(\Omega^{\text{ext}}\backslash B)]^{3n}$.
\end{definition}

\begin{remark}
	Obviously, in the same way as above, we can extend a Hölder continuous functions $p$ on 
	$\mathrm{cl}(D\backslash \Omega)$ to the whole domain $D$ for all $\Omega\in \mathcal{O}_{k,\alpha}$.
\end{remark}

With this definition of convergence, we can prove the compactness of $\mathcal{G}$. We show
that the metric spaces $(\mathcal{O}_{k,\alpha},d_H) ,\ 
(\mathcal{P}_{k,\alpha}, \lVert \cdot \rVert_{C^{k,\alpha^\prime}})$, and 
$(\mathcal{E}_{k,\alpha}, \lVert \cdot \rVert_{C^{k,\alpha^\prime}})$ are each compact, where $0 < \alpha^\prime < \alpha < 1$,
and that $\mathcal{G}$ is a closed subset of $\mathcal{O}_{k,\alpha}\times \mathcal{P}_{k,\alpha}\times \mathcal{E}_{k,\alpha}$.

\begin{lemma}
	\label{lemma:shape-compactness}
	The space of admissible shapes $\mathcal{O}_{k,\alpha}(\Omega_0,\Omega^{\text{ext}})$
	equipped with the \linebreak Hausdorff distance $d_H$ is a compact metric space.
\end{lemma}

\begin{proof}
	We prove that $\mathcal{O}_{k,\alpha}$ is sequentially compact. First, we show
	that 
	the space of design variables $U_{k,\alpha}^{\text{ad}}(\Omega^{\text{ext}})$ is compact.
	Then, the compactness of $\mathcal{O}_{k,\alpha}(\Omega_0,\Omega^{\text{ext}})$ follows out of
	it. By definition, $U_{k,\alpha}^{\text{ad}}(\Omega^{\text{ext}})$ is a bounded subspace of
	$C^{k,\alpha}(\Omega^{\text{ext}})$ and thus precompact in $C^{k,\alpha^\prime}(\Omega^{\text{ext}})$ for any 
	$0 < \alpha^\prime < \alpha$ (see Lemma \ref{lemma:pre-compactness}).
	Hence, for every sequence $(\psi_n)_{n\in \mathbb{N}}\subset
	U_{k,\alpha}^{\text{ad}}(\Omega^{\text{ext}})$ there exists a subsequence that converges in $C^{k,\alpha^\prime}(\Omega^{\text{ext}})$.
	For the compactness of $U_{k,\alpha}^{\text{ad}}(\Omega^{\text{ext}})$ it remains to show
	that the limit of $(\psi_n)_{n\in \mathbb{N}}$ lies in $U_{k,\alpha}^{\text{ad}}(\Omega^{\text{ext}})$.
	Since $U_{k,\alpha}^{\text{ad}}(\Omega^{\text{ext}})$ is precompact in the Banach space $C^{k,\alpha^\prime}(\Omega^{\text{ext}})$,
	that sequence has a subsequence
	$(\psi_{n_l})_{l\in\mathbb{N}}$ with $\psi_{n_l} \to \psi$ in $\lVert \cdot
	\rVert_{C^{k,\alpha^\prime}}$ for some $\psi\in C^{k,\alpha^\prime}(\Omega^{\text{ext}})$. First, we note that since
	$\lVert \psi_{n_l}
	\rVert_{C^{k,\alpha}}\leq K$ we have for any $\gamma \in
	\mathbb{N}^3$ with $0 \leq \lvert\gamma\rvert \leq k$
	\begin{equation*}
		\left\lvert \frac{\partial^{\lvert\gamma\rvert}\psi_{n_l}(x)}{\partial^{\gamma} x} -\frac{ \partial^{\lvert\gamma\rvert}\psi_{n_l}(y)}
		{\partial^{\gamma} x}\right\rvert \leq 
		\left(K - \max_{\lvert \gamma \rvert = k}
		\left\lVert \frac{\partial^{\lvert \gamma \rvert}\psi_{n_l}}
		{\partial x^\gamma} \right\rVert_{\infty} \right)
		\lvert x - y \rvert^\alpha
	\end{equation*}
	and
	\begin{equation*}
		\begin{split}
			&\max_{\lvert \gamma \rvert = k}\left[ \frac{\partial^{\lvert \gamma \rvert}\psi_{n_l}}
			{\partial x^\gamma} \right]_{0,\alpha} \leq
			K - \max_{\lvert \gamma \rvert \leq k}\left\lVert 
			\frac{\partial^{\lvert \gamma \rvert} \psi_{n_l}}{\partial x^\gamma} \right\rVert_\infty,\\
			&\max_{\lvert \gamma \rvert \leq k}\left\lVert 
			\frac{\partial^{\lvert \gamma \rvert} \psi_{n_l}}{\partial x^\gamma} \right\rVert_\infty
			\longrightarrow \max_{\lvert \gamma \rvert \leq k}\left\lVert 
			\frac{\partial^{\lvert \gamma \rvert} \psi}{\partial x^\gamma} \right\rVert_\infty
			\leq K.
		\end{split}
	\end{equation*}
	We apply these estimates to show that $\psi \in C^{k,\alpha}$ and $\lVert \psi \rVert_{C^{k,\alpha}}\leq K$:
	\begin{equation*}
		\begin{split}
			\left\lvert \frac{\partial^{\lvert\gamma\rvert}\psi(x)}{\partial^{\gamma} x} - 
			\frac{\partial^{\lvert\gamma\rvert}\psi(y)}
			{\partial^{\gamma} x}\right\rvert \leq
			\ &\left\lvert \frac{\partial^{\lvert\gamma\rvert}\psi(x)}{\partial^{\gamma} x} - 
			\frac{\partial^{\lvert\gamma\rvert}\psi_{n_l}(x)}
			{\partial^{\gamma} x}\right\rvert +
			\left\lvert \frac{\partial^{\lvert\gamma\rvert}\psi_{n_l}(x)}{\partial^{\gamma} x} - 
			\frac{\partial^{\lvert\gamma\rvert}\psi_{n_l}(y)}
			{\partial^{\gamma} x}\right\rvert \\
			+ \ &\left\lvert \frac{\partial^{\lvert\gamma\rvert}\psi_{n_l}(y)}{\partial^{\gamma} x} - 
			\frac{\partial^{\lvert\gamma\rvert}\psi(y)}
			{\partial^{\gamma} x}\right\rvert \\
			\leq \ &\left\lvert \frac{\partial^{\lvert\gamma\rvert}\psi(x)}{\partial^{\gamma} x} - 
			\frac{\partial^{\lvert\gamma\rvert}\psi_{n_l}(x)}
			{\partial^\gamma x}\right\rvert \\
			+ \
			&\left( K - \max_{\lvert \gamma \rvert \leq k}\left\lVert 
			\frac{\partial^{\lvert \gamma \rvert} \psi_{n_l}}{\partial x^\gamma} \right\rVert_\infty \right)
			\lvert x - y \rvert^\alpha \\
			+ \ &\left\lvert \frac{\partial^{\lvert\gamma\rvert}\psi_{n_l}(y)}{\partial^{\gamma} x} 
			- \frac{\partial^{\lvert\gamma\rvert}\psi(y)}{\partial^\gamma x}\right\rvert
		\end{split}
	\end{equation*}
	For the second term we used the Hölder continuity of $\psi_{n_l}$. The first
	and third term converge to zero since $\lVert \psi_{n_l} \rVert_{C^{k,\alpha}}
	\to \lVert \psi \rVert_{C^{k,\alpha}}$. Overall we get
	\begin{equation*}
		\begin{split}
			\left\lvert \frac{\partial^{\lvert\gamma\rvert}\psi(x)}{\partial^{\gamma} x} 
			- \frac{\partial^{\lvert\gamma\rvert}\psi(y)}
			{\partial^{\gamma} x}\right\rvert \leq
			&\left\lvert \frac{\partial^{\lvert\gamma\rvert}\psi(x)}{\partial^{\gamma} x} 
			- \frac{\partial^{\lvert\gamma\rvert}\psi_{n_l}(x)}
			{\partial^\gamma x}\right\rvert \\
			\ +
			&\left( K - \max_{\lvert \gamma \rvert \leq k}\left\lVert 
			\frac{\partial^{\lvert \gamma \rvert} \psi_{n_l}}{\partial x^\gamma} \right\rVert_\infty \right)
			\lvert x - y \rvert^\alpha \\
			+ \ &\left\lvert \frac{\partial^{\lvert\gamma\rvert}\psi_{n_l}(y)}{\partial^{\gamma} x} 
			- \frac{\partial^{\lvert\gamma\rvert}\psi(y)}{\partial^\gamma x}\right\rvert \\
			\longrightarrow \ &\left(K - \max_{\lvert \gamma \rvert \leq k}\left\lVert 
			\frac{\partial^{\lvert \gamma \rvert} \psi}{\partial x^\gamma} \right\rVert_\infty\right)
			\lvert x - y \rvert^\alpha.
		\end{split}
	\end{equation*}
	This gives $\psi \in C^{k,\alpha}(\Omega^{\text{ext}})$ and $\lVert \psi
	\rVert_{C^{k,\alpha}}\leq K$. Further, by the same arguments we can show that the sequence of inverse 
	$(\psi_{n_l}^{-1})_{l\in\mathbb{N}}$ converges to some function $\tilde{\psi} \in C^{k,\alpha}(\Omega^{\text{ext}})$ with respect to 
	$\lVert \cdot \rVert_{C^{k,\alpha^\prime}}$ (we can always find a subsequence) and with $\lVert \tilde{\psi} \rVert_{C^{k,\alpha}(\Omega^{\text{ext}})}
	\leq K$. 
	Now, it is straightforward to show that any bounded subset of $C^{k,\alpha}(\Omega^{\text{ext}})$ is uniformly equicontinuous, which implies that 
	$\tilde{\psi} = \psi^{-1}$ and thus $\psi \in U_{k,\alpha}^{\text{ad}}(\Omega^{\text{ext}})$.
	Therefore, $U_{k,\alpha}^{\text{ad}}(\Omega^{\text{ext}})$ is closed and with that compact. 
	
	We use the compactness of $U_{k,\alpha}^{\text{ad}}(\Omega^{\text{ext}})$ 
	with respect to $\lVert
	\cdot \rVert_{C^{k,\alpha^\prime}(\Omega^{\text{ext}})}$ to show the compactness of
	$(\mathcal{O}_{k,\alpha}, d_\mathcal{O})$. Consider a sequence $(\Omega_n)_{n\in \mathbb{N}} \subset
	\mathcal{O}_{k,\alpha}$. By the definition of $\mathcal{O}_{k,\alpha}$,
	there is a corresponding sequence $(\psi_n)_{n\in \mathbb{N}} \subset
	U^{\text{ad}}_{k,\alpha}(\Omega^{\text{ext}})$ with $\psi_n(\Omega_0) =
	\Omega_n$ for all $n\in \mathbb{N}$.
	Since $U^{\text{ad}}_{k,\alpha}(\Omega^{\text{ext}})$ is compact,
	there exists a subsequence $(\psi_{n_l})_{l\in \mathbb{N}}$ that converge
	to some $\psi \in U^{\text{ad}}_{k,\alpha}(\Omega^{\text{ext}})$ in
	$\lVert \cdot \rVert_{C^{k,\alpha^\prime}}$. We show
	that the corresponding subsequence of shapes $(\Omega_{n_l}) =
	\psi_{n_l}(\Omega_0)$ converge to $\Omega = \psi(\Omega_0)$ by using the
	convergence of $\psi_{n_l} \to \psi$ in $\lVert \cdot \rVert_{C^{k,\alpha\prime}}$:
	
	\begin{equation*}
		\begin{split}
			d_\mathcal{O}(\Omega_{n_l}, \Omega) = \ &\max\{\underset{x\in
				\Omega_{n_l}}{\sup} \ \underset{y\in \Omega}{\inf}\lvert x - y \rvert,\,
			\underset{y\in
				\Omega}{\sup} \ \underset{x\in \Omega_{n_l}}{\inf}\lvert x - y \rvert\}\\
			= \ &\max\{\underset{x\in
				\Omega_0}{\sup} \ \underset{y\in \Omega_0}{\inf}\lvert \psi_{n_l}(x) - \psi(y) \rvert,\,
			\underset{y\in
				\Omega_0}{\sup} \ \underset{x\in \Omega_0}{\inf}\lvert \psi_{n_l}(x) - \psi(y) \rvert\} \\
			\leq \ &\max\{\underset{x\in
				\Omega_0}{\sup}\lvert \psi_{n_l}(x) - \psi(x) \rvert,\,
			\underset{y\in
				\Omega_0}{\sup} \lvert \psi_{n_l}(y) - \psi(y) \rvert\} \\
			\underset{l\to \infty}{\longrightarrow} \ &0.
		\end{split}
	\end{equation*}
	Hence, each sequence in $\mathcal{O}_{k,\alpha}$ has a convergent
	subsequence that converge in $\mathcal{O}_{k,\alpha}$ w.r.t.
	the Hausdorff distance. Therefore, $(\mathcal{O}_{k,\alpha}, d_\mathcal{O})$ is
	sequentially compact. 
	
\end{proof}

\begin{lemma}
	\label{lemma:pot-compactness}
	Let $0 < \alpha^\prime < \alpha \leq 1$ and $k \geq 2$. Then, the solution space $\mathcal{P}_{k,\alpha}^{\mathrm{ext}}$,
	i.e., the space consisting of extensions from Definition \ref{def:extensions} to the solutions $\phi_\Omega \in
	\mathcal{P}_{k,\alpha}$, is compact in $C^{2,\alpha^\prime}(D)$.
\end{lemma}

\begin{proof}
	First, Lemma \ref{lemma:pot-eq} gives that $\mathcal{P}_{k,\alpha}$ consists of Hölder continuous functions. 
	We denote the extension from $\phi_\Omega  \in \mathcal{P}_{k,\alpha}$
	on $D$
	with $\phi^{\text{ext}}_\Omega$. With Lemma \ref{lemma:pot-uni} and (\ref{est:ext-bound}) 
	the following estimate holds
	\begin{equation*}
		\lVert \phi^{\text{ext}}_\Omega  \rVert_{C^{2,\alpha}(D)} \leq
		C\lVert \phi_{\Omega} \rVert_{C^{2,\alpha}(D\backslash\Omega)} \leq
		CK,
	\end{equation*}
	where $K$ is uniform in $\mathcal{O}_{k,\alpha}$. In \cite{BittGottsch} it is shown that
	the constant $C$ can also be chosen uniformly with respect to $\mathcal{O}_{k,\alpha}$ which
	yields an uniform bound for $\phi^{\text{ext}}_\Omega $. Hence, 
	$\mathcal{P}^\text{ext}_{k,\alpha}$ is a bounded subset of 
	$C^{k,\alpha^\prime}(D)$ and therefore precompact in 
	$C^{k,\alpha^\prime}(D)$ (see Lemma \ref{lemma:pre-compactness}). 
	Since $C^{k,\alpha^\prime}(D)$ is a 
	Banach space, it remains to show that $\mathcal{P}^\text{ext}_{k,\alpha}$ is closed.
	For this, let $(\phi_{\Omega_n}^\text{ext})_{n\in \mathbb{N}}\subset 
	\mathcal{P}^\text{ext}_{k,\alpha}$ be a sequence that converge to some function $\phi$ with respect to
	$\lVert \cdot \rVert_{C^{k,\alpha}(\Omega^{\text{ext}})}$ and let
	$(\Omega_n)_{n\in \mathbb{N}}\subset \mathcal{O}_{k,\alpha}$ be the corresponding sequence of shapes. As $\mathcal{O}_{k,\alpha}$
	is compact, we can find a subsequence of shapes $(\Omega_{n_l})_{l\in\mathbb{N}}$ that
	converge against some $\Omega\in \mathcal{O}_{k,\alpha}$. Now, consider the 
	corresponding subsequence of solutions $(\phi_{\Omega_{n_l}}^\text{ext})_{l\in \mathbb{N}}$. 
	This subsequence also converge against $\phi$. Since the convergence is in $\lVert \cdot \rVert_{C^{2,\alpha}(\Omega^{\text{ext}})}$
	and we have additionally seen in the proof of Lemma \ref{lemma:shape-compactness} that $\phi \in C^{2,\alpha}(D)$, it follows that $\phi$
	is the extension to a solution $\phi_\Omega$ for (\ref{eq:pot-eq}) and therefore lies in 
	$\mathcal{P}^\text{ext}_{k,\alpha}$. 
	
\end{proof}

\begin{lemma}
	\label{lemma:el-compactness}
	Let $0 < \alpha^\prime < \alpha \leq 1$ and $k \in \mathbb{N}_0$. The solution space $\mathcal{E}_{k,\alpha}^\mathrm{ext}$,
	i.e., the space consisting of extensions from Definition \ref{def:extensions} to the solutions $u_\Omega \in
	\mathcal{E}_{k,\alpha}$, is compact in $C^{2,\alpha^\prime}(\Omega^{\text{ext}})$.
\end{lemma}
\begin{proof}
	The proof follows the exact same arguments as in Lemma \ref{lemma:pot-compactness} and 
	therefore is omitted. 
	
\end{proof}

\begin{lemma}
	\label{lemma:graph-compactness}
	Consider the multi physics shape optimization problem {\normalfont (\ref{eq:multiphysics-shape-eq})} with boundary regularity 
	$C^{k,\alpha},\, k \geq 2$ and $\alpha \in \, ]0,1]$. Then, the Graph $\mathcal{G}$ is compact
	with respect to the corresponding maximum product metric.
\end{lemma}

\begin{proof}
	By applying Lemma \ref{lemma:shape-compactness}, Lemma \ref{lemma:pot-compactness}, and 
	Lemma \ref{lemma:el-compactness} we have that 
	$\mathcal{O}_{k,\alpha}\times 
	\mathcal{P}_{k,\alpha}^{\mathrm{ext}}\times
	\mathcal{E}_{k,\alpha}^{\mathrm{ext}}$ is compact. Let 
	$(\Omega_n)_{n\in \mathbb{N}} \subset \mathcal{O}_{k,\alpha}$
	and $\Omega \in \mathcal{O}_{k,\alpha}$ 
	with $\Omega_n \to \Omega$ in 
	$d_H$.  Then, one can see that due to the compactness
	of $\mathcal{O}_{k,\alpha}\times 
	\mathcal{P}_{k,\alpha}^{\mathrm{ext}}\times
	\mathcal{E}_{k,\alpha}^{\mathrm{ext}}$ that $\phi_{\Omega_n}^{\text{ext}} \to  
	\phi^{\text{ext}}$ and $u_{\Omega_n}^{\text{ext}} \to  
	u^{\text{ext}}$ in $\lVert \cdot \rVert_{C^{2,\alpha^\prime}(\Omega^{\text{ext}})}$
	with $\phi^{\text{ext}}_{\lvert \Omega}$ solves
	(\ref{eq:pot-eq}) on $\Omega$ and 
	$u^{\text{ext}}_{\lvert \Omega}$ 
	solves (\ref{eq:el-eq}) on $\Omega$. Hence,
	$\mathcal{G}$ is a closed subspace of a compact
	metric space and therefore compact as well. 
	
\end{proof}

\begin{lemma}[Continuity of Local Cost Functionals \cite{GottschSchmitz}]
	\label{lemma:cost-funct-continuity}
	Let $\mathcal{F}_{\text{vol}},\, \mathcal{F}_{\text{sur}}
	\in C^0(\mathbb{R}^d)$ (with $d$ as in Definition \ref{def:cost-functional} with $r=3$) 
	and let $\mathcal{O}_{k,\alpha}$ only consists $C^0$-admissible 
	shapes. For $\Omega$ and $\boldsymbol{v} \in [C^k(\mathrm{cl}(\Omega))]^{3n}$
	consider the volume integral $J_{\text{vol}}(\Omega,\boldsymbol{v})$ and
	the surface integral $J_{\text{sur}}(\Omega,\boldsymbol{v})$.	
	Let $\Omega_n\subset\mathcal{O}_{k,\alpha}$ with $\Omega_n
	\overset{\tilde{\mathcal{O}}}{\longrightarrow}\Omega$
	as $n\to \infty$ and let $(\boldsymbol{v}_n)_{n\in\mathbb{N}}\subset
	[C^k(\mathrm{cl}(\Omega_n))]^{3n}$ be a sequence with 
	$\boldsymbol{v}_n \rightsquigarrow \boldsymbol{v}$ as $n\to \infty$ for some 
	$\boldsymbol{v}\in [C^k(\mathrm{cl}(\Omega_n))]^{3n}$. Then,
	\begin{itemize}
		\item[(i)] $J_{\text{vol}}(\Omega_n,\boldsymbol{v}_n) \longrightarrow
		J_{\text{vol}}(\Omega,\boldsymbol{v})$ as $n\to \infty$.
		\item[(ii)] If the set $\mathcal{O}_{k,\alpha}$ only consists of $C^1$-admissible shapes one obtains \\
		$J_{\text{sur}}(\Omega_n,\boldsymbol{v}_n)$  $\longrightarrow
		J_{\text{sur}}(\Omega,\boldsymbol{v})$ as $n\to \infty$ as well.
	\end{itemize}
\end{lemma}

\begin{proof}
	(i) First, we apply the characteristic function on the volume integral and obtain
	\begin{equation*}
		J_{\text{vol}}(\Omega_n,\boldsymbol{v}_n) := \int_{\Omega^\text{ext}} 
		\chi_{\Omega_n}\cdot\mathcal{F}_{\text{vol}}\hspace{.1em} 
		(x,\boldsymbol{v}^\text{ext}_n,\nabla \boldsymbol{v}^\text{ext}_n,\dots, 
		\nabla^k \boldsymbol{v}^\text{ext}_n)\,dx.
	\end{equation*}
	Because of $\mathcal{F}_{\text{vol}} \in C^0(\mathbb{R^d})$ and $\boldsymbol{v}_n
	\rightsquigarrow \boldsymbol{v}$ as $n\to \infty$ there exist a constant $C>0$ such
	that $\lvert\chi_{\Omega_n}\cdot\mathcal{F}_{\text{vol}}\hspace{.1em} 
	(x,\boldsymbol{v}^\text{ext}_n,\nabla \boldsymbol{v}^\text{ext}_n,
	\dots, \nabla^k \boldsymbol{v}^\text{ext}_n)\rvert
	\leq C$ is valid for all $n\in \mathbb{N}$ almost everywhere in $\Omega_{\text{ext}}$.
	Moreover, $\Omega_n \overset{\tilde{\mathcal{O}}}{\longrightarrow} \Omega$ 
	and $\boldsymbol{v}_n^{\text{ext}} \to \boldsymbol{v}^{\text{ext}}$ in 
	$[C_0^k(\Omega^\text{ext})]^{3n}$ ensure the existence of  
	\begin{equation*}
		\begin{split}
			&\lim_{n\to\infty}\chi_{\Omega_n}\cdot\mathcal{F}_{\text{vol}}\hspace{.1em} 
			(x,\boldsymbol{v}^\text{ext}_n,
			\nabla \boldsymbol{v}^\text{ext}_n,\dots, \nabla^k \boldsymbol{v}^\text{ext}_n) \\
			&= \chi_{\Omega}\cdot\mathcal{F}_{\text{vol}}\hspace{.1em} 
			(x,\boldsymbol{v}^\text{ext},
			\nabla \boldsymbol{v}^\text{ext},\dots, \nabla^k \boldsymbol{v}^\text{ext}),
		\end{split}
	\end{equation*}
	for all $x\in\Omega^{\text{ext}}$. Therefore, we can apply Lebesgue`s dominated convergence theorem:
	\begin{equation*}
	\begin{split}
		\lim_{n\to\infty} J_{\text{vol}}(\Omega_n,\boldsymbol{v}_n) &= 
		\lim_{n\to\infty} \int_{\Omega^\text{ext}} 
		\chi_{\Omega_n}\cdot\mathcal{F}_{\text{vol}}\hspace{.1em} 
		(x,\boldsymbol{v}^\text{ext}_n,\nabla \boldsymbol{v}^\text{ext}_n,\dots, 
		\nabla^k \boldsymbol{v}^\text{ext}_n)\,dx \\
		&= \ \int_{\Omega^\text{ext}}  \lim_{n\to\infty} 
		\chi_{\Omega_n}\cdot\mathcal{F}_{\text{vol}}\hspace{.1em} 
		(x,\boldsymbol{v}^\text{ext}_n,\nabla \boldsymbol{v}^\text{ext}_n,\dots, 
		\nabla^k \boldsymbol{v}^\text{ext}_n)\,dx \\
		&= \ \int_{\Omega^\text{ext}}
		\chi_{\Omega}\cdot\mathcal{F}_{\text{vol}}\hspace{.1em} 
		(x,\boldsymbol{v}^\text{ext},\nabla \boldsymbol{v}^\text{ext},\dots, 
		\nabla^k \boldsymbol{v}^\text{ext})\,dx \\
		&= \ J_{\text{vol}}(\Omega,\boldsymbol{v})
	\end{split}
	\end{equation*}
	(ii) The second assertion can be analogously proven as in \cite{BittGottsch} and we only state
	the main ideas here.
	
	First, we note that every shape $\Omega \in \mathcal{O}_{k,\alpha}$, by its definition, can be
	considered as two-dimensional submanifold and therefore is locally embeddable into $\mathbb{R}^2$.
	Let $A_n^i \subset \partial\Omega$, $1\leq i \leq m$ with
	$\cup_{i=1}^m A_i = \partial\Omega$ and 
	$A_i \cap A_j = \emptyset$ for $i\neq j$. We can find in $i$ and $n$ uniformly bounded 
	chart mappings $h_n^i:A_n^i \to \tilde{A}_i$ with $\tilde{A}_i\subset \mathbb{R}^2$. 
	We use them to straighten
	the boundary of $\Omega_n$ to obtain a volume integral, e.g.,
	\begin{equation*}
	\begin{split}
		J_{\text{sur}}(\Omega_n,\boldsymbol{v}_n) &= \ \int_{\partial \Omega_n} 
		\mathcal{F}_{\text{sur}}\hspace{.1em}
		(x,\boldsymbol{v}_n,\nabla \boldsymbol{v}_n,\dots, \nabla^k \boldsymbol{v}_n)\,dA \\
		&= \ \sum_{i=1}^{m} \int_{A_n^i} 
		\mathcal{F}_{\text{sur}}\hspace{.1em}
		(x,\boldsymbol{v}_n,\nabla \boldsymbol{v}_n,\dots, \nabla^k \boldsymbol{v}_n)\,dA \\
		&\hspace{-1.5cm}= \ \sum_{i=1}^{m} \int_{\tilde{A}_i} 
		\mathcal{F}_{\text{sur}}\hspace{.1em}
		\left(h_n^i(s),\boldsymbol{v}_n(h_n^i(s)),
		\nabla \boldsymbol{v}_n(h_n^i(s)),\dots, 
		\nabla^k \boldsymbol{v}_n(h_n^i(s))\right)\sqrt{g^{h_n^i}(s)}\,ds,
	\end{split}
	\end{equation*}
	where $ds$ indicates the Lebesgue measure on $\tilde{A}_i$
	and corresponding Gram determinants $g^{h_n^i}$. Due to the fact that the chart mappings 
	$h_n^i$ are uniform bounded and since $\tilde{A}_i$ is independent of $n$, one can see that,
	similarly to (i), we can apply Lebesgue's Theorem and the assertion is proven.
	
\end{proof}

\begin{remark}
	The continuity assumption of Lemma \ref{lemma:cost-funct-continuity} ensures
	the existence of an integrable majorant for $\mathcal{F}_{\text{vol}}$ and  
	$\mathcal{F}_{\text{sur}}$. Example (\ref{functional:energy})
	does not fulfil this assumption. However, (\ref{coefficient:share-stress})
	is integrable on compact sets and one can easily find an integrable majorant
	by applying the uniform bound of Lemma \ref{lemma:pot-uni}.
\end{remark}

\begin{theorem}
	Assuming $k \geq 2$, the multi physics shape optimization problem { \normalfont (\ref{eq:multiphysics-shape-eq})} possesses
	at least one Pareto optimal solution $(\Omega^*, \phi_{\Omega^*}, u_{\Omega^*})
	\in \mathcal{G}$ and covers all nondominated points in $\mathcal{Y}$, i.e.,
	$\mathcal{Y}_N = \mathrm{cl}({\mathcal{Y}})_N$.
\end{theorem}
\begin{proof}
	Lemma \ref{lemma:graph-compactness} provides the compactness of $\mathcal{G}$ and 
	Lemma \ref{lemma:cost-funct-continuity} the continuity of the local cost functionals.
	Then, Theorem \ref{theorem:shape-solution-existence} gives the existence of an optimal shape 
	and the closeness of the set of optimal shapes.
	
\end{proof}

\section{\label{section:scalarization}Scalarization and Multi Physics Optimization}

Scalarizing is the traditional approach to solving a multicriteria optimization problem.
This includes formulating a single objective optimization problem that is related to 
the original Pareto optimality problem by means of a real-valued scalarizing function
typically being a function of the objective function, auxiliary scalar or vector variables,
and/or scalar or vector parameters. Additionally scalarization techniques sometimes
further restrict the feasible set of the problem with new variables or/and restriction
functions.
In this section we investigate the stability of the parameter-dependent optimal shapes to
different types of scalarization techniques with underlying design problem
(\ref{eq:multiphysics-shape-eq}).

First, let us define the scalarization methods we consider. This involves a certain
class of real-valued functions $S_\theta : \mathbb{R}^l\to \mathbb{R}$, referred to 
as scalaraization function that possibly depends on a parameter $\theta$ which lies in a parameter
space $\Theta$. The scalarization problem is given by
\begin{equation}
\label{eq:scal-eq}
\begin{split}
\begin{array}{ll}
\min S_\theta\left(\mathbf{J}(\Omega,u_{\Omega},\phi_{\Omega})\right)\\
\text{\hspace{.4em} subject to } (\Omega,u_{\Omega},\phi_{\Omega}) \in 
\mathcal{G}_\theta,
\end{array}
\end{split}
\end{equation}
where $\mathcal{G}_\theta \subseteq \mathcal{G}$.
For the sake of notational convenience, we sometimes identify an element
$(\Omega, u_\Omega, \phi_\Omega)\in \mathcal{G}_\theta$ only by its distinct shape $\Omega$.
If we assume that $\mathcal{G}_\theta$ is closed and
the scalarization $S_\theta(\mathbf{J})$ is lower semicontinuous on
$\mathcal{G}_\theta \times \{ \theta \}$,
then, by the results of Section \ref{section:existence}, (\ref{eq:scal-eq}) obviously
has an optimal solution for $\theta\in\Theta$.
For a fixed $\theta\in\Theta$ we shall denote the space of all optimal shapes 
to an achievement function problem with 
$\zeta_\theta=\argmin_{\Omega\in \mathcal{G}_\theta}S_\theta\left(\mathbf{J}(\Omega,u_{\Omega},\phi_{\Omega})\right)$. We assume that 
$\Theta \subset \mathbb{R}^l$ is closed and equip the space
$\mathcal{Z}:= \{\zeta_\theta \, :\, \theta\in\Theta\}$ with the Hausdorff distance,
which in this setting defines, due to the closeness of the optimal shapes sets, a metric
(see Lemma \ref{lemma:chi-close} and Corollary \ref{cor:m-space}).

In the following, we gather some definitions and assertion from chapter 4 of 
\cite{bank83}. 
We define the optimal set mapping $\chi: \Theta \rightrightarrows
\mathcal{Z}$, the optimal value mapping 
$\tau:\Theta \longrightarrow \mathbb{R}$, and the graph mapping $G:\Theta
\rightrightarrows 2^\mathcal{G}$ which maps a parameter
$\theta\in\Theta$ to the corresponding set of optimal shapes 
$\zeta_\theta$, the corresponding optimal value $\min_{\Omega\in\mathcal{G}_\theta} S_\theta(\mathbf{J})$,
and the corresponding graph $\mathcal{G}_\theta$, respectively.
With these definitions in hand, we can describe the stability of the optimal shapes for a wide range
of scalarization methods. First, we state a lemma that shows that $(\mathcal{Z}, d_H)$ is indeed a
metric space.

\begin{lemma}
\label{lemma:chi-close}
    The optimal set mapping $\chi$ is closed if $\tau$ is upper semicontinuous and
    $S_\theta(\mathbf{J})$ is lower semicontinuous on $\mathcal{G}\times \{ \theta \}$.
\end{lemma}

\begin{corollary}
\label{cor:m-space}
    If the scalarization function $S_\theta$ is lower semicontinuous on $\mathbb{R}^l \times \{ \theta \}$
    and uniform continuous on $\{ r \} \times \Theta$, for $r\in \mathbb{R}^l$, then the Hausdorff 
    distance $d_H$ defines a metric on $\mathcal{Z}$. 
\end{corollary}
\begin{proof}
    Due to the continuity of $\mathbf{J}$ (see Lemma \ref{lemma:cost-funct-continuity})
    and the uniform  continuity of $S_\theta$ on $\{ r \} \times \Theta$, the optimal value maping
    $\tau$ is upper semicontinuous and therefore, by Lemma \ref{lemma:chi-close}, the optimal set 
    mapping $\chi$ is closed. Since $d_H$ defines a metric on $F(\mathcal{G})$ (the set of all
	closed subsets of $\mathcal{G}$), $(\mathcal{Z}, d_H)$ defines a metric space.
	
\end{proof}

Since the sclarization solution is not necessarly unique, we need some sort of 
continuity property of point-to-set mappings in order to discuss the stability of sets of optimal shapes. 
The literatur describes
serveral definitions which vary in the statement.
We investigate the stabilty according to Hausdorff and Berge (for Berge see \cite{bank83}) which, in this setting,
are equivalent.

\begin{definition}[Upper semicontiniuity according to Hausdorff]
	Let $(\Theta, d_\Theta)$ and $(X, d_x)$ be metric spaces. A point-to-set mapping
	of $\Theta$ into $X$ is a function $\Gamma$ that assigns a subset $\Gamma(\theta)$
	of $X$ to each element $\theta\in\Theta$. This function is called upper semicontinuous in 
	$\theta^*$ if for each sequence $(\theta_n)_{n\in \mathbb{N}} \subseteq \Theta$ with 
	$\theta_n \longrightarrow \theta^*$, for $n \to \infty$, we have 
	\begin{equation}
		\label{eq:optimal-shapes-convergens}
		\sup_{x\in T(\theta_n)} 
		\inf_{x^\prime\in T(\theta^*)} 
		d_X(x, x^\prime) \longrightarrow 0.
	\end{equation}
	$\Gamma$ is called upper semicontinuous if $\Gamma$ is upper semicontinuous in each 
	$\theta\in \Theta$. 
	For this type of continuity we simply write u.s.c.-H.
\end{definition}

The next Theorem states stability conditions for scalarization function problems.

\begin{theorem}[\cite{bank83}]
	\label{theorem:argmin-conv}
	Assume that $G$ is u.s.c.-H at $\theta^*$ and $G(\theta^*)$ is compact. Further, 
	let $\tau$ be upper semicontinuous at $\theta^*$ and $S_{\theta^*}$ lower semicontinuous
	on $G(\theta^*)\times \{ \theta^* \}$. Thenm the optimal set mapping $\chi$ is 
	u.s.c.-H at $\theta^*$.
\end{theorem}

The following two corollaries demonstrate continuity properties of shapes under change of preferences for two commonly used 
scalarization techniques. In particular, the results apply to the shape optimization problem introduced in Section 
\ref{section:problem}. 

\begin{corollary}[Weighted Sum Method]
	Consider cost functionals \linebreak $\mathbf{J} =  (J_1,\dots, J_l)$ and let $\Theta\subset\mathbb{R}^l$ be a
	closed subset. Then, the weighted sum scalaraization method, which is given by
	\begin{equation*}
		\begin{split}
		\begin{array}{ll}
		\min \sum_{i=1}^l \theta_i J_i((\Omega,u_{\Omega},\phi_{\Omega}))\\
		\text{\hspace{.4em} subject to } (\Omega,u_{\Omega},\phi_{\Omega}) \in 
		\mathcal{G},
		\end{array}
		\end{split}
	\end{equation*}
	fulfills all conditions of Theorem \ref{theorem:argmin-conv} due to the compactness of $\mathcal{G}$ 
	(see Lemma \ref{lemma:graph-compactness}) and the continuity of $\mathbf{J}$ (see Lemma \ref{lemma:cost-funct-continuity}).
\end{corollary}

\begin{corollary}[$\epsilon$-Constraint Method]
	Let $\mathbf{J} = (J_1,\dots, J_l)$ be cost functionals. We optimize 
	cost functional $J_j$ and constrain the other functionals by $J_i \leq \epsilon_i \in \mathbb{R}$, for
	$1 \leq i \leq n$ and $i\neq j$. If each $\epsilon_i$ converges monotonically decreasing to some
	$\epsilon_i^*$, then the $\epsilon$-Constraint Method
	\begin{equation*}
		\begin{split}
		\begin{array}{ll}
		\min J_j((\Omega,u_{\Omega},\phi_{\Omega}))\\
		\text{\hspace{.4em} subject to } J_i \leq \epsilon_i,
		\end{array}
		\end{split}
	\end{equation*}
	fulfills all conditions of Theorem \ref{theorem:argmin-conv}.
\end{corollary}
\begin{proof}
	Let $\epsilon = (\epsilon_1, \dots, \epsilon_l)$ and $\mathcal{G}_\epsilon = \{ \Omega\in\mathcal{G} \ : \ 
	J_i(\Omega) \leq \epsilon_i,\, i\neq j\}$. The u.s.c.-H of $G$ is given due to the continuity of $\mathbf{J}$.
	The continuity of $J_j$, the u.s.c.-H of $G$, and the fact that $\mathcal{G}_\epsilon \subseteq \mathcal{G}_{\epsilon^\prime}$ for all
	$\epsilon^*\leq\epsilon^\prime \leq \epsilon$
	gives that $\tau(\epsilon)$ converge continuously against $\tau(\epsilon^*)$ for $\epsilon \searrow \epsilon^*$. Hence, the optimal sets 
	$\chi(\epsilon)$ converge against $\chi(\epsilon^*)$ for $\epsilon \searrow \epsilon^*$ in the sense of u.s.c.-H.
	
\end{proof}
\begin{remark}
	Whenever the scalarized problem \eqref{eq:scal-eq} possesses a unique solution $\zeta_\theta=\{\Omega_\theta\}$ for all $\theta$ in 
	some neighborhood of $\theta^*\in\Theta$, then $\Omega_{\theta_n}\longrightarrow \Omega_{\theta^*}$ in the Hausdorff distance 
	(for subsets in $\mathbb{R}^d$) if $\theta_n\to\theta^*$.  
\end{remark}

\section{\label{section:outlook} Conclusions}
In this work, we extended the well known framework of the existence of optimal solutions in shape optimization \cite{Chenais75,hasmaek03} 
to a multi criteria setting. We formulated conditions for the existence and completeness of Pareto optimal points.  Multiple criteria in 
design are often related to simulations that include different domains of physics. We presented a coupled fluid-dynamic and mechanical system, 
which is motivated by gas turbine design and fits to the given framework. The objective functions in this case are given by fluid losses and 
mechanical durability expressed by the probability of failure under low cycle fatigue. Both objectives require classical solutions to the 
underlying partial differential equations and therefore can only be formulated on sufficiently regular shapes such that elliptic regularity 
theory applies \cite{GilbTrud,Agm59,Agm64}. We presented a formulation of the family of admissible shapes that implied the existence of such 
classical solutions and thereby provided a non-trivial example for the general framework. 

In multi criteria optimization \cite{ehrgott2005multicriteria}, the Pareto front contains points which are optimal with respect to different 
preferences of a decision maker. An interesting point is, if a small variation of the preference also leads to a small variation in the design. 
This question however is ill-posed, if Pareto optimal solutions need not to be unique. We therefore presented a study, where we went over to the sets 
of Pareto optimal shapes for a given preference and studied the variation of these sets in the Hausdorff-metric. In this setting, certain 
continuity properties in the preference parameters were derived.

It will be of interest to develop multi-criteria shape optimization also from an algorithmic standpoint, using the theory of shape derivatives and 
gradient based optimization; see, e.g., \cite{doganay19,bolten20} for some first steps in that direction.  For a rigorous analysis of numerical 
schemes of shape optimization, it will be of interest, if (a) the optima of the discretized problem are close to the optima of the continuous problem 
and, if (b) the same holds for shape gradients for non-optimal solutions, as e.g., used in multi-criteria descent algorithms. In particular, this 
should be true for the objective values of discretized and continuous solutions, respectively. Potentially, iso-geometric finite elements 
\cite{cottrell2009isogeometric,Fusseder2015isogeometric,Wall2008isogeometric} could be an useful numerical tool to not spoil the domain regularity 
that is built into our framework by the need of $C^{k,\alpha}$-classical solutions needed for the evaluation of the objectives in multi-criteria 
shape optimization problems like the one presented here.

\section*{Acknowledgements} 
   The authors gratefully acknowledge interesting discussions with 
   Laura Bittner and Kathrin Klamroth. This work has been funded by the federal 
   ministry of research (BMBF) via the GIVEN project, grant no. 05M18PXA.

\appendix

\section{\label{app}Appendix}
\begin{definition}[Hölder continuity]\label{def:holder_continuity}
	Let $U \subseteq \mathbb{R}^d$ be open. A function $f:U \to \mathbb{R}$ is called
	Hölder continuous if there exist non-negative real constants $C > 0,\, \alpha \in \, ]0,1]$,
	such that
	\begin{equation*}
		\lvert f(x) - f(y) \rvert \leq C \lvert x - y \rvert^\alpha, \quad \text{for all } x,y \in U.
	\end{equation*} 
	With $C^{k,\alpha}(U)$ we denote the space consisting of function on $U$ having
	Hölder continuous derivatives up to order $k \in \mathbb{N}_{0}$ with exponent $\alpha$. 
	If the function $f$ and its derivatives up to order $k$ are bounded on the 
	closure of $U$, then we can assign the norm
	\begin{equation*}
		\lVert f \rVert_{C^{k,\alpha}\left(U\right)} = 
		\max_{\lvert \gamma \rvert \leq k}\left\lVert 
		\frac{\partial^{\lvert \gamma \rvert} f}{\partial x^\gamma} \right\rVert_\infty
		+ \max_{\lvert \gamma \rvert = k}\left[ \frac{\partial^{\lvert \gamma \rvert}f}
		{\partial x^\gamma} \right]_{0,\alpha},
	\end{equation*}
	where
	\begin{equation*}
		\begin{split}
			&\lVert f \rVert_{\infty} = \sup_{x\in U} \lvert f(x) \rvert  \\
			&[f]_{0,\alpha} = \sup_{\substack{x,y \in U \\ x \neq y}} 
			\frac{\lvert f(x) - f(y) \rvert}{\lvert x - y \rvert^\alpha}.
		\end{split}
	\end{equation*}
	We note that $(C^{k,\alpha}(U), \lVert \cdot \rVert_{C^{k,\alpha}(U)})$ is a Banach space. 
	A vector field $f:U \to \mathbb{R}^n$ is Hölder continuous with exponent 
	$\alpha$ if each component of $f$ is Hölder continuous with exponent $\alpha$.
\end{definition} 

\begin{lemma}[\cite{GilbTrud}, Lemma 6.37]
	\label{lemma:extension}
	Let $\Omega$ be a $C^{k,\alpha}$ domain in $\mathbb{R}^d$ with $k\geq 1, \, \alpha \in \, ]0,1]$ and let
	$\Omega^\prime$ be an open set containing $\mathrm{cl}({\Omega})$. Suppose 
	$u\in C^{k,\alpha}(\mathrm{cl}\left({\Omega}\right))$. Then, there exist a function 
	$w\in C_0^{k,\alpha}(\Omega^\prime)$ such that $w = u$ and
	\begin{equation}
	\label{est:ext-bound}
	\lVert w \rVert_{C^{k,\alpha}(\Omega^\prime)} \leq
	C\lVert u \rVert_{C^{k,\alpha}(\Omega)},
	\end{equation}
	where $C=C(k,\Omega, \Omega^\prime)$.
\end{lemma}

\begin{lemma}[Korn's Second Inequality, \cite{braess97}]
	\label{lemma:Korn-2nd}
	Let $\Omega \subset \mathbb{R}^3$ be an open bounded set with piecewise smooth
	boundary. In addition, suppose $\Gamma_0 \subset \partial\Omega$ has positive
	two dimensional measure. Then, there exist a positive number $c^\prime =
	c^\prime(\Omega, \Gamma_0)$ such that
	\begin{equation*}
		\lVert \epsilon(u_{\Omega}) \rVert^2_{L_2(\Omega)} = \int_\Omega \mathrm{tr}(\epsilon(v)^2)\, dx 
		\geq c^\prime \lVert v \rVert_{H^1(\Omega)}, \quad \text{for all } v\in H^1_\Gamma(\Omega).
	\end{equation*}
	Here, $H^1_\Gamma(\Omega)$ is the closure of $\{ v \in [C^\infty(\Omega)]^3 \,
	:\, v(x) = 0 \text{ for } x\in \Gamma_0 \}$ with respect to the 
	$\lVert \cdot \rVert_{H^1(\Omega)}$-norm.
\end{lemma}

\begin{lemma}[\cite{GilbTrud}, Lemma 6.36]
	\label{lemma:pre-compactness}
	Let $\Omega$ be a $C^{k,\alpha}$ domain in $\mathbb{R}^d$ with $k\geq 1, \, \alpha\in \, ]0,1]$ and let
	$S$ be a bounded set in $C^{k,\alpha}(\mathrm{cl}\left({\Omega}\right))$. Then, $S$ is precompact in 
	$C^{j,\beta}(\mathrm{cl}\left({\Omega}\right))$ if $j+\beta \leq k + \alpha$.
\end{lemma}

\bibliographystyle{spmpsci_unsrt}
\bibliography{references}

\end{document}